\documentclass[11pt,reqno]{amsart}
\usepackage{latexsym,amssymb,amsmath,multicol,rotating,
lscape}
\font\sm=msbm7 scaled \magstep 2
\font\germ=eufm10 scaled 1200

\newtheorem{thm}{Theorem}[section]

\newtheorem{prop}[thm]{Proposition}

\newcommand{\thmref}[1]{Theorem~\ref{#1}}
\newcommand{\lemref}[1]{Lemma~\ref{#1}}
\newcommand{\rmkref}[1]{Remark~\ref{#1}}
\newcommand{\propref}[1]{Proposition~\ref{#1}}
\newcommand{\corref}[1]{Corollary~\ref{#1}}
\newcommand{\C}{{\mathbb C}}
\theoremstyle{remark}

\begin{document}

\title[octonary quadratic forms]
{On the number of representations by certain octonary  quadratic forms with coefficients 1, 2, 3, 4 and 6}
 
\author{B. Ramakrishnan, Brundaban Sahu and Anup Kumar Singh}
\address[B. Ramakrishnan and Anup Kumar Singh]{Harish-Chandra Research Institute, HBNI, 
       Chhatnag Road, Jhunsi,
     Allahabad -     211 019,
   India.}
\address[Brundaban Sahu]
{School of Mathematical Sciences, National Institute of Science 
Education and Research, HBNI, Bhubaneswar, Via- Jatni, Khurda, Odisha - 752 050,
India.}

\email[B. Ramakrishnan]{ramki@hri.res.in}
\email[Brundaban Sahu]{brundaban.sahu@niser.ac.in}
\email[Anup Kumar Singh]{anupsingh@hri.res.in}

\date{\today}

\begin{abstract}
In this paper, we find formulas for the  number of representations of certain diagonal octonary quadratic forms with coefficients $1,2,3,4$ and $6$. We obtain these 
formulas by constructing explicit bases of the space of modular forms of weight $4$ on $\Gamma_0(48)$ with character. 
\end{abstract}

\maketitle

\section{Introduction}

Let ${\mathbb N}, {\mathbb N}_0$ and ${\mathbb Z}$ denote the set of positive integers, non-negative integers and integers respectively. 
For $a_1, \ldots, a_8 \in {\mathbb N}$ and $n\in {\mathbb N}_0$, we define 
$$
N(a_1, \ldots, a_8;n) := {\rm card}\left\{(x_1,\ldots, x_8)\in {\mathbb Z}^8 \vert n = a_1 x_1^2 + \cdots + a_8 x_8^2\right\}.
$$
Note that $N(a_1, \ldots, a_8;0) =1$. Without loss of generality we may assume that 
$$
a_1\le a_2\le \cdots \le a_8 {\rm ~and~} \gcd(a_1, \ldots, a_8) =1.
$$
Formulae for $N(a_1,\ldots, a_8;n)$ for the octonary quadratic forms 
\begin{equation}\label{1236}
\sum_{r=1}^i x_r^2 + 2\sum_{r=i+1}^{i+j} x_r^2 + 3\sum_{r=i+j+1}^{i+j+k} x_r^2+ 6\sum_{r=i+j+k+1}^{i+j+k+l} x_r^2  
\end{equation}
for all the partitions $i+j+k+l =8$, $i,j,k,l\ge 0$ appeared in the literature. When all of them ($i,j,k,l$) are even (there are 26 cases), it was obtained by 
several authors (see \cite{{a-a-w}, {a-k1}, {a-w}, {jacobi}, {kokluce1}, {kokluce2}, {ijnt1}}). When all of them are odd (there are 10 cases), it was obtained in \cite{a-k2}. Recently, in \cite{a-k-new}, the authors considered the rest of the cases (mixed parity) and this completed all the cases for the coefficients $1, 2, 3$ or $6$. In \cite{kokluce3} the author considers one octonary quadratic form with coefficients $2,3,6$ and $12$. A few cases of the coefficients $1,2,3,6$ 
are also considered in this paper. It is to be noted that various methods were used to obtain these formulas such as elementary evaluations, 
using the method of convolution sums of the divisor functions and the theory of modular forms. However, in most of the cases techniques involving  modular forms were used. 

Formulae for $N(a_1, \ldots, a_8;n)$  for the octonary quadratic forms 
\begin{equation}\label{124}
\sum_{r=1}^i x_r^2 + 2\sum_{r=i+1}^{i+j} x_r^2 + 4\sum_{r=i+j+1}^{i+j+k} x_r^2, 
\end{equation}
where $i+j+k =8$, $(i,j,k) \in$ $\{(7,0,1), (1,0,7), (6,0,2), (2,0,6), (5,2,1), (1,2,5), (5,0,3), (3,0,5), \\
(4,2,2), (2,2,4), (2,4,2), (4,0,4), (3,4,1), (3,2,3), (1,6,1), (1,4,3)\}$  were obtained in \cite{{a-a-w}, {a-a-w1}}. There are a total of $38$ cases with coefficients 
$1,2,4$ and out of which $10$ of them are obtainable from the previous works. The remaining $12$ cases are $\{(1,1,6), (1,3,4), (1,5,2), (2,1,5), (2,3,3), (2,5,1),\\
(3,1,4), (3,3,2),  (4,1,3), (4,3,1), (5,1,2),  (6,1,1) \}$. 

In this paper, we first complete the $12$ cases remaining for the coefficients $1,2,4$ and then extend the above works to find formulae for the diagonal octonary  quadratic forms with coefficients $1, 2, 3, 4$ and $6$ by using the theory of modular forms. More precisely, we consider the following octonary quadratic forms (with coefficients $1,2,3,4,6$):
\begin{equation}\label{12346}
\sum_{r=1}^i x_r^2 + 2\sum_{r=i+1}^{i+j} x_r^2 + 3\sum_{r=i+j+1}^{i+j+k} x_r^2+ 4\sum_{r=i+j+k+1}^{i+j+k+l} x_r^2  + 6\sum_{r=i+j+k+l+1}^{i+j+k+l+m} x_r^2,  
\end{equation}
where $i+j+k+l+m =8$ and find their number of representations $N(a_1,\ldots, a_8;n)$. The corresponding theta series will be a modular form of weight $4$ on $\Gamma_1(48)$. The case $l=0$ is the case with coefficients $1,2,3,6$, which was done earlier (as mentioned above) and so we take $l\not= 0$. When $m=0$, there are 84 cases (with $k,l\not =0$) which is given in Table 1. (Here we have not included the cases $k=0$, which corresponds to the coefficients $1,2,4$  and $l=0$, which corresponds to the coefficients $1,2,3$.)  When $m\not= 0$, there are 210 cases which is given in Table 2, except for  the following $7$ cases ($(0,0,0,1,7), (0,0,0,2,6), (0,0,0,3,5)$, $(0,0,0,4,4), (0,0,0,5,3), (0,0,0,6,2), (0,0,0,7,1)$), which  can be obtained from the earlier results with coefficients $2,3$ only. So, here we consider only the remaining 203 cases. In both the tables, we indicate the corresponding modular forms spaces. We also remark that the $21$ cases $\{(0,4,0,2,2), (0,2,0,4,2), (0,2,0,2,4), (0,5,0,1,2), (0,1,0,3,4), 
\linebreak
(0,3,0,3,2), (0,1,0,5,2), (0,3,0,1,4), (0,1,0,1,6), (0,3,0,2,3),(0,1,0,4,3), (0,3,0,4,1), (0,5,0,2,1), \\
(0,1,0,2,5), (0,1,0,6,1), (0,4,0,1,3), (0,2,0,3,3), (0,6,0,1,1), (0,2,0,1,5), (0,4,0,3,1), (0,2,0,5,1)\}$ can also be obtained from the results arising from the coefficients $1,2,3$. However, these cases are 
also kept in the respective tables. Since we consider  different bases for the modular forms spaces, our formulas in these cases are also 
different from the previous formulas obtained by Alaca and Kesicio\u{g}lu \cite{{a-k-new}, {a-k1}}. 


\section{Statement of the Results}

Let us first consider the 12 remaining cases of the quadratic forms given by \eqref{124} with $i+j+k=8$ as mentioned in the introduction. We denote the number of representations of an integer $n$ by these quadratic forms as $N(1^i,2^j,4^k;n)$. Then, the following theorem gives the 
representation numbers $N(1^i,2^j,4^k;n)$ for the 12 cases $(i,j,k) \in \{(1,1,6), (1,3,4), (1,5,2), (2,1,5), (2,3,3), (2,5,1), (3,1,4), \\
 (3,3,2),  (4,1,3), (4,3,1), (5,1,2),  (6,1,1) \}$. 


\begin{thm}\label{1}
Let $n\in {\mathbb N}$. Then 
\begin{equation*}
\begin{split}
{\rm (i)} ~N(1^1,2^1,4^6;n) &= \frac{2}{11} \sigma_{3;\chi_8, {\bf 1}}(n) +\frac{2}{11} \sigma_{3;{\bf 1},\chi_8}(n/2) + \frac{6}{11}  a_{4,8,\chi_8;1}(n) + \frac{14}{11} a_{4,8,\chi_8;2}(n) \\
& \quad  + \frac{48}{11}  a_{4,8,\chi_8;1}(n/2) - \frac{28}{11} a_{4,8,\chi_8;2}(n/2),\\
{\rm (ii)} ~N(1^1,2^3,4^4;n) &= \frac{4}{11} \sigma_{3;\chi_8, {\bf 1}}(n) +\frac{2}{11} \sigma_{3;{\bf 1},\chi_8}(n/2) + \frac{1}{11}  a_{4,8,\chi_8;1}(n) + \frac{17}{11} a_{4,8,\chi_8;2}(n) \\
& \quad  + \frac{48}{11}  a_{4,8,\chi_8;1}(n/2) +\frac{16}{11} a_{4,8,\chi_8;2}(n/2),\\
{\rm (iii)} ~N(1^1,2^5,4^2;n) &= \frac{8}{11} \sigma_{3;\chi_8, {\bf 1}}(n) +\frac{2}{11} \sigma_{3;{\bf 1},\chi_8}(n/2) + \frac{2}{11}  a_{4,8,\chi_8;1}(n) + \frac{12}{11} a_{4,8,\chi_8;2}(n) \\
& \quad  + \frac{48}{11}  a_{4,8,\chi_8;1}(n/2) +\frac{16}{11} a_{4,8,\chi_8;2}(n/2),\\
{\rm (iv)} ~N(1^2,2^1,4^5;n) &= \frac{4}{11} \sigma_{3;\chi_8, {\bf 1}}(n) +\frac{2}{11} \sigma_{3;{\bf 1},\chi_8}(n/2) + \frac{12}{11}  a_{4,8,\chi_8;1}(n) + \frac{28}{11} a_{4,8,\chi_8;2}(n) \\
& \quad  + \frac{92}{11}  a_{4,8,\chi_8;1}(n/2) - \frac{28}{11} a_{4,8,\chi_8;2}(n/2),\\
\end{split}
\end{equation*}
\begin{equation*}
\begin{split}
{\rm (v)} ~N(1^2,2^3,4^3;n) &= \frac{8}{11} \sigma_{3;\chi_8, {\bf 1}}(n) +\frac{2}{11} \sigma_{3;{\bf 1},\chi_8}(n/2) + \frac{2}{11}  a_{4,8,\chi_8;1}(n) + \frac{34}{11} a_{4,8,\chi_8;2}(n) \\
& \quad  + \frac{92}{11}  a_{4,8,\chi_8;1}(n/2) +\frac{16}{11} a_{4,8,\chi_8;2}(n/2),\\
{\rm (vi)} ~N(1^2,2^5,4^1;n) &= \frac{16}{11} \sigma_{3;\chi_8, {\bf 1}}(n) +\frac{2}{11} \sigma_{3;{\bf 1},\chi_8}(n/2) + \frac{4}{11}  a_{4,8,\chi_8;1}(n) + \frac{24}{11} a_{4,8,\chi_8;2}(n) \\
& \quad  + \frac{48}{11}  a_{4,8,\chi_8;1}(n/2) +\frac{16}{11} a_{4,8,\chi_8;2}(n/2),\\
{\rm (vii)} ~N(1^3,2^1,4^4;n) &= \frac{8}{11} \sigma_{3;\chi_8, {\bf 1}}(n) +\frac{2}{11} \sigma_{3;{\bf 1},\chi_8}(n/2) + \frac{13}{11}  a_{4,8,\chi_8;1}(n) + \frac{45}{11} a_{4,8,\chi_8;2}(n) \\
& \quad  + \frac{136}{11}  a_{4,8,\chi_8;1}(n/2) +\frac{16}{11} a_{4,8,\chi_8;2}(n/2),\\
{\rm (viii)} ~N(1^3,2^3,4^2;n) &= \frac{16}{11} \sigma_{3;\chi_8, {\bf 1}}(n) +\frac{2}{11} \sigma_{3;{\bf 1},\chi_8}(n/2) + \frac{4}{11}  a_{4,8,\chi_8;1}(n) + \frac{46}{11} a_{4,8,\chi_8;2}(n) \\
& \quad  + \frac{136}{11}  a_{4,8,\chi_8;1}(n/2) +\frac{16}{11} a_{4,8,\chi_8;2}(n/2),\\
{\rm (ix)} ~N(1^4,2^1,4^3;n) &= \frac{16}{11} \sigma_{3;\chi_8, {\bf 1}}(n) +\frac{2}{11} \sigma_{3;{\bf 1},\chi_8}(n/2) + \frac{4}{11}  a_{4,8,\chi_8;1}(n) + \frac{68}{11} a_{4,8,\chi_8;2}(n) \\
& \quad  + \frac{180}{11}  a_{4,8,\chi_8;1}(n/2) +\frac{104}{11} a_{4,8,\chi_8;2}(n/2),\\
{\rm (x)} ~N(1^4,2^3,4^1;n) &= \frac{32}{11} \sigma_{3;\chi_8, {\bf 1}}(n) +\frac{2}{11} \sigma_{3;{\bf 1},\chi_8}(n/2) + \frac{8}{11}  a_{4,8,\chi_8;1}(n) + \frac{48}{11} a_{4,8,\chi_8;2}(n) \\
& \quad  + \frac{136}{11}  a_{4,8,\chi_8;1}(n/2) +\frac{16}{11} a_{4,8,\chi_8;2}(n/2),\\
{\rm (xi)} ~N(1^5,2^1,4^2;n) &= \frac{32}{11} \sigma_{3;\chi_8, {\bf 1}}(n) +\frac{2}{11} \sigma_{3;{\bf 1},\chi_8}(n/2) - \frac{14}{11}  a_{4,8,\chi_8;1}(n) + \frac{92}{11} a_{4,8,\chi_8;2}(n) \\
& \quad  + \frac{224}{11}  a_{4,8,\chi_8;1}(n/2) +\frac{192}{11} a_{4,8,\chi_8;2}(n/2),\\
{\rm (xii)} ~N(1^6,2^1,4^1;n) &= \frac{64}{11} \sigma_{3;\chi_8, {\bf 1}}(n) +\frac{2}{11} \sigma_{3;{\bf 1},\chi_8}(n/2) - \frac{28}{11}  a_{4,8,\chi_8;1}(n) + \frac{96}{11} a_{4,8,\chi_8;2}(n) \\
& \quad  + \frac{224}{11}  a_{4,8,\chi_8;1}(n/2) +\frac{192}{11} a_{4,8,\chi_8;2}(n/2).\\
\end{split}
\end{equation*}
The terms appearing on the right-hand side of the above formulas are defined in \S 4.1.
\end{thm}

\smallskip

Next we shall state the formulae for the quadratic forms with coefficients $1,2,3,4$ given in Table 1. We state them as four statements in the theorem, each statement corresponds to 
the four modular forms spaces that appear in the Table 1. 

\begin{thm}\label{2}
Let $n\in {\mathbb N}$ and $i,j,k,l$ be non-negative integers such  that $i+j+k+l =8$. \\
{\rm (i)} For each entry $(i,j,k,l)$ in Table {\rm 1} corresponding to the space $M_4(\Gamma_0(48))$, we have 
\begin{equation}
N(1^i,2^j,3^k,4^l;n) = \sum_{\alpha=1}^{30} a_\alpha A_\alpha(n),
\end{equation}
where $A_\alpha(n)$ are the Fourier coefficients of the basis elements $f_\alpha$ defined in \S {\rm 4.2} and the values of the constants $a_\alpha$ are 
given in Table {\rm 3}. \\
{\rm (ii)} For each entry $(i,j,k,l)$ in Table $1$ corresponding to the space $M_4(\Gamma_0(48),\chi_8)$, we have 
\begin{equation}
N(1^i,2^j,3^k,4^l;n) = \sum_{\alpha=1}^{28} b_\alpha B_\alpha(n),
\end{equation}
where $B_\alpha(n)$ are the Fourier coefficients of the basis elements $g_\alpha$ defined in \S {\rm 4.3} and the values of the constants $b_\alpha$ are 
given in Table {\rm 4}. \\
{\rm (iii)} For each entry $(i,j,k,l)$ in Table $1$ corresponding to the space $M_4(\Gamma_0(48),\chi_{12})$, we have 
\begin{equation}
N(1^i,2^j,3^k,4^l;n) = \sum_{\alpha=1}^{30} c_\alpha C_\alpha(n),
\end{equation}
where $C_\alpha(n)$ are the Fourier coefficients of the basis elements $h_\alpha$ defined in \S {\rm 4.4} and the values of the constants $c_\alpha$ are 
given in Table {\rm 5}. \\
{\rm (iv)} For each entry $(i,j,k,l)$ in Table $1$ corresponding to the space $M_4(\Gamma_0(48),\chi_{24})$, we have 
\begin{equation}
N(1^i,2^j,3^k,4^l;n) = \sum_{\alpha=1}^{28} d_\alpha D_\alpha(n),
\end{equation}
where $D_\alpha(n)$ are the Fourier coefficients of the basis elements $F_\alpha$ defined in \S {\rm 4.5} and the values of the constants $d_\alpha$ are 
given in Table {\rm 6}. 
\end{thm}

\bigskip

In the following theorem we list the formulas for the octonary quadratic forms with coefficients 
$1,2,3,4,6$ corresponding to Table 2. The proof is similar to the proof of \thmref{2} and so we omit the details. 

\smallskip

\begin{thm}\label{3}
Let $n\in {\mathbb N}$ and $i,j,k,l,m$ be non-negative integers such  that $i+j+k+l+m =8$. \\
{\rm (i)} For each entry $(i,j,k,l,m)$ in Table {\rm 2} corresponding to the space $M_4(\Gamma_0(48))$, we have 
\begin{equation}
N(1^i,2^j,3^k,4^l,6^m;n) = \sum_{\alpha=1}^{30} a'_\alpha A_\alpha(n),
\end{equation}
where $A_\alpha(n)$ are the Fourier coefficients of the basis elements $f_\alpha$ defined in \S {\rm 4.2} and the values of the constants $a'_\alpha$ are 
given in Table {\rm 7}. \\
{\rm (ii)} For each entry $(i,j,k,l,m)$ in Table $2$ corresponding to the space $M_4(\Gamma_0(48),\chi_8)$, we have 
\begin{equation}
N(1^i,2^j,3^k,4^l,6^m;n) = \sum_{\alpha=1}^{28} b'_\alpha B_\alpha(n),
\end{equation}
where $B_\alpha(n)$ are the Fourier coefficients of the basis elements $g_\alpha$ defined in \S {\rm 4.3} and the values of the constants $b'_\alpha$ are 
given in Table {\rm 8}. \\
{\rm (iii)} For each entry $(i,j,k,l,m)$ in Table $2$ corresponding to the space $M_4(\Gamma_0(48),\chi_{12})$, we have 
\begin{equation}
N(1^i,2^j,3^k,4^l,6^m;n) = \sum_{\alpha=1}^{30} c'_\alpha C_\alpha(n),
\end{equation}
where $C_\alpha(n)$ are the Fourier coefficients of the basis elements $h_\alpha$ defined in \S {\rm 4.4} and the values of the constants $c'_\alpha$ are 
given in Table {\rm 9}. \\
{\rm (iv)} For each entry $(i,j,k,l,m)$ in Table $2$ corresponding to the space $M_4(\Gamma_0(48),\chi_{24})$, we have 
\begin{equation}
N(1^i,2^j,3^k,4^l,6^m;n) = \sum_{\alpha=1}^{28} d'_\alpha D_\alpha(n),
\end{equation}
where $D_\alpha(n)$ are the Fourier coefficients of the basis elements $F_\alpha$ defined in \S {\rm 4.5} and the values of the constants $d'_\alpha$ are 
given in Table {\rm 10}. 
\end{thm}

\smallskip

\section{Preliminaries}

As we use the theory of modular forms, we shall first present some preliminary facts on modular forms. For $k\in \frac{1}{2}{\mathbb Z}$,
let $M_k(\Gamma_0(N),\chi)$ denote the space of modular forms of weight $k$ for the congruence subgroup $\Gamma_0(N)$ with character $\chi$ and $S_k(\Gamma_0(N), \chi)$ be the subspace of cusp forms of weight $k$ for $\Gamma_0(N)$ with character $\chi$. We assume $4\vert N$ when $k$ is not an integer and in that case, the character $\chi$ which is a Dirichlet character modulo $N$, is an even character. When $\chi$ is the trivial (principal) character 
modulo $N$, we shall denote the spaces by $M_k(\Gamma_0(N))$ and $S_k(\Gamma_0(N))$ respectively. Further, when $k\ge 4$ is an integer and $N=1$, we shall denote these vector spaces by $M_k$ and $S_k$ respectively. 

For an integer $k \ge 4,$ let $E_k$ denote the normalized Eisenstein series of weight $k$ in $M_k$ given by 
$$
E_k(z) = 1 - \frac{2k}{B_k}\sum_{n\ge 1} \sigma_{k-1}(n) q^n,
$$
where $q=e^{2 i\pi z}$, $z\in {\mathcal H}$, the complex upper half-plane, $\sigma_r(n)$ is the sum of the $r$th powers of the positive divisors of $n$, and $B_k$ is the $k$-th Bernoulli number defined by $\displaystyle{\frac{x}{e^x-1} = \sum_{m=0}^\infty \frac{B_m}{m!} x^m}$.

The classical theta function which is fundamental to the theory of modular forms of half-integral weight is defined by 
\begin{equation}\label{theta}
\Theta(z) = \sum_{n\in {\mathbb Z}} q^{n^2},
\end{equation}
and is a modular form in the space $M_{1/2}(\Gamma_0(4))$. Another function which is mainly used in our work is the Dedekind eta function $\eta(z)$ and 
it is given by 
\begin{equation}\label{eta}
\eta(z)=q^{1/24} \prod_{n\ge1}(1-q^n).
\end{equation}

An eta-quotient is a finite product of integer powers of $\eta(z)$ and we denote it as follows: 
\begin{equation}\label{eta-q} 
\prod_{i=1}^s \eta^{r_i}(d_i z) := d_1^{r_1} d_2^{r_2} \cdots d_s^{r_s},
\end{equation}
where $d_i$'s are positive integers and $r_i$'s are non-zero integers.

In the following we shall present some facts about modular forms of integral and half-integral weights, which we shall be using in our proof. 

\noindent {\bf Fact I}.\\
{We give a fact about certain duplication of modular forms, which follow from the two results: Chapter 3, Propostion 17 of  \cite{koblitz} and 
 Proposition 1.3 of \cite{shimura}. 
If $f$ is a modular form in $M_k(\Gamma_0(N), \chi)$, then for a positive integer $d$, the function $f(dz)$ is a modular form in $M_k(\Gamma_0(dN), \chi)$, 
if $k$ is an integer and it belongs to the space $M_k(\Gamma_0(dN), \chi \chi_d)$, if $k$ is a half-integer.
Here, we have used the notation $\chi_m := \left(\frac{m}{\cdot}\right)$, the Kronecker symbol, where $m$ is a non-zero integer, which is a 
character modulo $|m|$. 
}

\noindent {\bf Fact II}.\\
{For  positive integers $r$, $r_1$, $r_2$,  $d_1$, $d_2$, we have 
\begin{equation}\label{theta1}
\Theta^{r}(d_1z) \in 
\begin{cases} M_{r/2}(\Gamma_0(4d_1), \chi_{d_1}) & {\rm ~if~} r {\rm ~is~odd}, \\ 
                           M_{r/2}(\Gamma_0(4d_1), \chi_{-4}) & {\rm ~if~} r \equiv 2\pmod{4},\\ 
                           M_{r/2}(\Gamma_0(4d_1)) & {\rm ~if~} r \equiv 0\pmod{4}.
\end{cases}   
\end{equation}

\begin{equation}\label{theta2}
\Theta^{r_1}(d_1z)\cdot \Theta^{r_2}(d_2z)  \in 
\begin{cases} 
M_{\frac{r_1+r_2}{2}}(\Gamma_0(4 [d_1, d_2]), \chi_{(-d_1d_2)}) & \!\!{\rm ~if~} r_1 r_2 {\rm ~is~odd~},  r_1+r_2\equiv 2\pmod{4}, \\
M_{\frac{r_1+r_2}{2}}(\Gamma_0(4 [d_1, d_2]), \chi_{(d_1d_2)}) &  \!\!{\rm ~if~} r_1 r_2 {\rm ~is~odd~}, r_1+r_2 \equiv 0\pmod{4}.
\end{cases}   
\end{equation}

In order to get the above fact we use the following properties. By Fact I, if $f\in M_k(\Gamma_0(N),\chi)$, then $f(dz)$ belongs to $M_k(\Gamma_0(dN), \chi')$, where $\chi'=\chi$ if $k$ is an integer and $\chi' = \chi \chi_d$, if $k$ is a half-integer. Next, if $f_i\in M_{k_i}(\Gamma_0(4N_i), \chi_i)$, $i=1,2$ are two modular forms of weight $k_i$ (where $k_1$ and $k_2$ are integers or half-integers). Then, 
it follows that the product $f_1 f_2$ is a modular form in $M_{k_1+k_2}(\Gamma_0(4[N_1,N_2]), \psi)$, where $\psi$ is a character modulo 
$4[N_1,N_2]$. If both the weights $k_1$ and $k_2$ are integers, then the resulting form is of weight $k_1+k_2$ (which is an integer) and so 
$\psi(-1) = (-1)^{k_1+k_2} = \chi_1(-1) \chi_2(-1)$, implying $\psi = \chi_1 \chi_2$. If $k_1+k_2$ is half-integer (say) with $k_1$ integer and 
$k_2$ half-integer, then $\chi_1(-1) = (-1)^{k_1}$ and $\chi_2$ is an even character modulo $4N_2$. Since the resulting form is of half-integral weight 
we must have $\psi$ an even character. Therefore, $\psi = \chi_1\chi_2$ if $k_1$ is even and $\psi = \chi_1\chi_2\chi_{-4}$ if $k_1$ is odd. 
If both $k_1$ and $k_2$ are half-integer, then the resulting form is of integer weight $k_1+k_2$. Since  $\chi_1\chi_2$  
is an even character,  when $k_1+k_2$ is an odd integer, the character of the space should be $\chi_1\chi_2\chi_{-4}$ instead of $\chi_1\chi_2$. 
Using these facts, we get the special cases as mentioned in \eqref{theta1} and \eqref{theta2}. For details we refer to \cite[Chap. 4, Proposition 3]{koblitz} and \cite[Proposition 1.3]{shimura}. 
}

\noindent {\bf Fact III}.\\
It is a fact that the vector space $M_k(\Gamma_1(N))$ is decomposed into modular forms space with character as follows. For this fact we refer to 
\cite[Proposition 28, p. 137]{koblitz}. 
\begin{equation}
M_k(\Gamma_1(N)) = \oplus_{\chi}M_k(\Gamma_0(N), \chi),
\end{equation}
where the direct sum varies over all Dirichlet characters modulo $N$ if the weight $k$ is a  positive integer and varies over all even 
Dirichlet characters modulo $N$, $4\vert N$, if the weight $k$ is half-integer. Further, if $k$ is an integer, one has $M_k(\Gamma_0(N),\chi) = \{0\}$, if $\chi(-1) \not= (-1)^k$. We also have the following decomposition of the space into subspaces of Eisenstein series and cusp forms:
\begin{equation}
M_k(\Gamma_0(N),\chi) = {\mathcal E}_k(\Gamma_0(N),\chi) \oplus S_k(\Gamma_0(N), \chi),
\end{equation}
where ${\mathcal E}_k(\Gamma_0(N),\chi)$ is the space generated by the Eisenstein series of weight $k$ on $\Gamma_0(N)$ with character $\chi$. 

\noindent {\bf Fact IV}.\\
For this fact we refer to the works of Atkin-Lehner and Li \cite{{a-l}, {li}}. 
By the Atkin-Lehner theory of newforms, the space $S_k(\Gamma_0(N),\chi)$ can be decomposed into the space of newforms and oldforms:
\begin{equation}
S_k(\Gamma_0(N),\chi) = S_k^{new}(\Gamma_0(N),\chi) \oplus S_k^{old}(\Gamma_0(N),\chi),
\end{equation}
where the above is an orthogonal direct sum (with respect to the Petersson scalar product) and  
\begin{equation}
S_k^{old}(\Gamma_0(N), \chi) = \bigoplus_{r\vert N, r<N\atop{rd\vert N}}S_k^{new}(\Gamma_0(r),\chi)\vert B(d).
\end{equation}
In the above, $S_k^{new}(\Gamma_0(N),\chi)$ is the space of newforms and $S_k^{old}(\Gamma_0(N),\chi)$ is the space of oldforms and the operator 
$B(d)$ is given by $f(z) \mapsto f(dz)$.

\noindent {\bf Fact V}:\\
Suppose that $\chi$ and $\psi$ are primitive Dirichlet characters with conductors $M$ and $N$, respectively. For a positive integer $k$, let 
\begin{equation}\label{eisenstein}
E_{k,\chi,\psi}(z) :=  c_0 + \sum_{n\ge 1}\left(\sum_{d\vert n} \psi(d) \cdot \chi(n/d) d^{k-1}\right) q^n,
\end{equation}
where 
$$
c_0 = \begin{cases}
0 &{\rm ~if~} M>1,\\
- \frac{B_{k,\psi}}{2k} & {\rm ~if~} M=1,
\end{cases}
$$
and $B_{k,\psi}$ denotes generalized Bernoulli number with respect to the character $\psi$. 
Then, the Eisenstein series $E_{k,\chi,\psi}(z)$ belongs to the space $M_k(\Gamma_0(MN), \chi/\psi)$, provided $\chi(-1)\psi(-1) = (-1)^k$ 
and $MN\not=1$. When $\chi=\psi =1$ (i.e., when $M=N=1$) and $k\ge 4$, we have $E_{k,\chi,\psi}(z) = - \frac{B_k}{2k} E_k(z)$, where $E_k$ is the normalized Eisenstein series of weight $k$ as defined before. 
For more details we refer to \cite[Chapter 7]{miyake} and \cite[Section 5.3]{stein}.  

We give a notation to the inner sum in \eqref{eisenstein}:
\begin{equation}\label{divisor}
\sigma_{k-1;\chi,\psi}(n) := \sum_{d\vert n} \psi(d) \cdot \chi(n/d) d^{k-1}.
\end{equation}

\smallskip

For more details on the theory of modular forms of integral and half-integral weights, we refer to \cite{{a-l}, {koblitz}, {li}, {miyake}, {shimura}}.

\smallskip

\section{Proofs}

In this section, we shall give a proof of our results. As mentioned in the introduction, we shall be using the theory of modular forms. 
Using Fact II, it is easy to see that the theta series associated to the quadratic forms with coefficients $1,2,4$ belong to the space 
$M_4(\Gamma_0(16),\chi)$, where $\chi$ is a Dirichlet character modulo $16$ and the theta series associated to the quadratic forms with coefficients $1,2,3,4,6$ belong to the space 
$M_4(\Gamma_0(48),\psi)$, where $\psi$ is a Dirichlet character modulo $48$. Therefore, in order to get the required formulae for $N(a_1,\ldots, a_8;n)$, we need a basis for these spaces. 
We shall give explict bases for the following spaces of modular forms of weight $4$:\\
$$
M_4(\Gamma_0(16), \chi_8), M_4(\Gamma_0(48)), M_4(\Gamma_0(48),\chi_8), M_4(\Gamma_0(48),\chi_{12}), M_4(\Gamma_0(48),\chi_{24}).
$$
(We have used the $L$-functions and modular forms database \cite{lfmdb} to get some of the cusp forms of weight $4$.)
 
\subsection{A basis for $M_4(\Gamma_0(16), \chi_8)$ and  proof of \thmref{1}.}\label{4.1}
The vector space $M_4(\Gamma_0(16),\chi_8)$ has dimension $8$ and the cusp forms space $S_4(\Gamma_0(16),\chi_8)$ has dimension $4$. 
Moreover, we have  \\
$S_4^{new}(\Gamma_0(16),\chi_8) = \{0\}$ and $S_4^{new}(\Gamma_0(8),\chi_8)$ is 2-dimensional. Let ${\bf 1}$ denote the trivial character with 
conductor $1$. Then by Fact V, the Eisenstein series $E_{4,{\bf 1}, \chi_8}(z)$ and $E_{4,\chi_8, {\bf 1}}(z)$ span the space 
${\mathcal E}_4(\Gamma_0(8),\chi_8)$. They are given by 
\begin{equation}\label{eis:8chi8}
E_{4, {\bf 1},\chi_8}(z) ~=~ \frac{11}{2} + \sum_{n\ge 1} \sigma_{3;{\bf 1},\chi_8}(n) q^n, \quad E_{4,\chi_8, {\bf 1}}(z) ~=~  \sum_{n\ge 1} \sigma_{3;\chi_8, {\bf 1}}(n) q^n.
\end{equation}
The space $S_4^{new}(\Gamma_0(8),\chi_8)$ is spanned by the following two eta-quotients (from now onwards we will be using the notation given in \eqref{eta-q}):
\begin{equation}\label{8chi8}
f_{4,8,\chi_8;1}(z) = 1^{-2}  2^{11} 4^{-3}  8^2 
= \sum_{n\ge 1} a_{4,8,\chi_8;1}(n) q^n ,~~ 
f_{4,8,\chi_8;2}(z) = 1^{2}  2^{-3} 4^{11}  8^{-2}  
=  \sum_{n\ge 1} a_{4,8,\chi_8;2}(n) q^n. 
\end{equation}
In the following proposition, we shall give a basis of $M_4(\Gamma_0(16),\chi_8)$.

\begin{prop}\label{prop1}
A basis of $M_4(\Gamma_0(16),\chi_8)$ is given by 
\begin{equation}
\left\{E_{4,{\bf 1},\chi_8}(tz), E_{4,\chi_8, {\bf 1}}(tz), f_{4,8,\chi_8;1}(tz),  f_{4,8,\chi_8;2}(tz)\big\vert t\mid 2\right\}.
\end{equation}
\end{prop}

\noindent 
We are now ready to prove \thmref{1}.

In these cases, the quadratic forms have coefficients $1,2, 4$. Therefore, by Fact II, all of them belong to the space of modular forms 
of weight $4$ on $\Gamma_0(16)$. Since the power of theta function corresponding to the coefficient $2$ is odd in all these cases, the modular forms spaces 
will have charcter $\chi_8$. Therefore, expressing each of the theta series corresponding to these 12 quadratic forms as a linear combination of 
the basis elements given in \propref{prop1}. Now, by comparing the $n$-th Fourier coefficients, we obtain the formulas listed in \thmref{1}.

\smallskip

\subsection{A basis for $M_4(\Gamma_0(48))$ and proof of \thmref{2}(i).}

The vector space $M_4(\Gamma_0(48))$ has dimension $30$ and we have $\dim_{\mathbb C}{\mathcal E}_4(\Gamma_0(48)) = 12$ and 
$\dim_{\mathbb C}S_4(\Gamma_0(48)) = 18$. For $d=6,8,12,16$ and $24$, $S_4^{new}(\Gamma_0(d))$ is one-dimensional and 
$\dim_{\mathbb C}S_4^{new}(\Gamma_0(48)) =3$. Let us define some eta-quotients and use them to give an explicit basis for $S_4(\Gamma_0(48))$.
Let
\begin{eqnarray}
& f_{4,6}(z) = 1^{2}  2^{2} 3^{2}  6^2 := \displaystyle{\sum_{n\ge 1}} a_{4,6}(n) q^n, \quad 
 f_{4,8}(z) = 2^4  4^4 := \displaystyle{\sum_{n\ge 1}} a_{4,8}(n) q^n, \\
& f_{4,12}(z) =  
1^{-1} 2^2 3^3 4^3 6^2 12^{-1} - 1^3 2^2 3^{-1} 4^{-1} 6^2 12^3 := \displaystyle{\sum_{n\ge 1}} a_{4,12}(n) q^n, \\
& f_{4,16}(z) = 2^{-4} 4^{16} 8^{-4} := \displaystyle{\sum_{n\ge 1}} a_{4,16}(n) q^n, ~
f_{4,24}(z) =  1^{-4} 2^{11} 3^{-4} 4^{-3} 6^{11} 12^{-3} := \displaystyle{\sum_{n\ge 1}} a_{4,24}(n) q^n.
\end{eqnarray}

Let $\chi_{-4}$ be the primitive odd character modulo $4$. Then the following new Eisenstein series belongs to ${\mathcal E}_4(\Gamma_0(16))$:
\begin{equation}\label{eis16} 
E_{4,\chi_{-4},\chi_{-4}}(z) = \sum_{n\ge 1} \sigma_{3,\chi_{-4},\chi_{-4}}(n) q^n = \sum_{n \ge 1} \left(\frac{-4}{n}\right) \sigma_3(n) q^n.
\end{equation}
We use the following notation in the sequel. For a Dirichlet character $\chi$ and a function $f$ with Fourier expansion $f(z) =\sum_{n\ge 1} a(n) q^n$, we define the twisted function $f \otimes \chi (z)$ as follows. 
\begin{equation}\label{twist}
f\otimes \chi (z) = \sum_{n\ge 1} \chi(n) a(n)q^n.
\end{equation}

A basis for the space $M_4(\Gamma_0(48))$ is given in the following proposition. 

\begin{prop}\label{trivial}
A basis for the Eisenstein series space ${\mathcal E}_4(\Gamma_0(48))$ is given by 
\begin{equation}
\left\{E_{4,\chi_{-4},\chi_{-4}}(z),E_{4,\chi_{-4},\chi_{-4}}(3z), E_4(tz), t\vert 48 \right\}  
\end{equation}
and a basis for the space of cusp forms $S_4(\Gamma_0(48))$ is given by 
\begin{equation}
\begin{split}
\left\{f_{4,6}(t_1z), t_1\vert 8; f_{4,8}(t_2z), t_2\vert 6; f_{4,12}(t_3z), t_3\vert 4;  f_{4,16}(t_4z),t_4\vert 3; \right. &\\ 
\left.  f_{4,24}\otimes \chi_4(t_5z), t_5\vert 2; f_{4,6}\otimes \chi_{-4}(z), f_{4,12}\otimes \chi_{-4}(z), f_{4,24}\otimes \chi_{-4}(z),\right\}&\\
\end{split}
\end{equation}
Together they form a basis for $M_4(\Gamma_0(48))$. 
\end{prop}

For the sake of simplicity in the formulae, we list these basis elements as $\{f_\alpha(z)\vert 1\le\alpha\le30\}$, where  $f_1(z) = E_4(z)$, $ f_2(z) = E_4(2z) $, $ f_3(z) = E_4(3z)$,
$ f_4(z) = E_4(4z)$, $f_5(z) = E_4(6z)$, $ f_6(z) = E_4(8z)$, $ f_7(z) = E_4(12z)$, $ f_{8}(z) = E_4(16z)$, $ f_{9}(z) = E_4(24z)$, $f_{10}(z) = E_4(48z)$, $f_{11}(z) = E_{4,\chi_{-4},\chi_{-4}}(z)$, $f_{12}(z) = E_{4,\chi_{-4},\chi_{-4}}(3z)$,  $f_{13}(z) = f_{4,6}(z)$, $f_{14}(z) = f_{4,6}(2z)$, $f_{15}(z) = f_{4,6}(4z)$, $f_{16}(z) = f_{4,6}(8z)$,  $f_{17}(z) = f_{4,8}(z)$, $f_{18}(z) = 
f_{4,8}(2z)$, $f_{19}(z) = f_{4,8}(3z)$, $f_{20}(z) = f_{4,8}(6z)$, $f_{21}(z) = f_{4,12}(z)$, 
$f_{22}(z) = f_{4,12}(2z)$, $f_{23}(z) = f_{4,12}(4z)$,   $f_{24}(z) = f_{4,16}(z)$, $f_{25}(z) = f_{4,16}(3z)$,  $f_{26}(z) = f_{4,24}\otimes \chi_4 (z)$,  $f_{27}(z) = f_{4,24}\otimes \chi_4 (2z)$,
$f_{28}(z) = f_{4,6}\otimes \chi_{-4}(z)$, $f_{29}(z) =  f_{4,12}\otimes \chi_{-4}(z)$, 
$f_{30}(z) =  f_{4,24}\otimes \chi_{-4}(z)$. 

\noindent We also express the Fourier coefficients of the function $f_\alpha(z)$ as $\sum_{n\ge 1} A_\alpha(n) q^n$, $1\le \alpha\le 30$. 

We are now ready to prove the theorem. Noting that all the 10 cases corresponding to the trivial character space (in Table 1) have the property 
that the powers of the theta functions corresponding to the coefficients $2$ and $3$ are even. Therefore, the resulting functions belong to the 
space of modular forms of weight $4$ on $\Gamma_0(48)$ with trivial character (we use Fact II to prove this). So, we can express these theta 
functions as a linear combination of the basis given in \propref{trivial} as follows. 
\begin{equation}
\Theta^i(z) \Theta^j(2z)\Theta^k(3z)\Theta^l(4z) = \sum_{\alpha=1}^{30} a_\alpha f_\alpha(z),
\end{equation}
where $a_\alpha$'s some constants. Comparing the $n$-th Fourier coefficients on both the sides, we get 
\begin{equation*}
N(1^i,2^j,3^k,4^l;n) = \sum_{\alpha=1}^{30} a_\alpha A_\alpha(n).
\end{equation*}
Explicit values for the constants \!$a_\alpha$, \!$1\le\alpha\le 30$ corresponding to the 10 cases are given in Table \!3. 

\subsection{A basis for $M_4(\Gamma_0(48), \chi_8)$ and proof of \thmref{2}(ii)} 

The space $M_4(\Gamma_0(48), \chi_8)$ is $28$ dimensional and the cusp forms space has dimension $20$. 
For the space of Eisenstein series we use the basis elements of ${\mathcal E}_4(\Gamma_0(8),\chi_8)$ given in \eqref{eis:8chi8}. 
For the space of cusp forms, there are no newforms and the oldforms classes are $S_4^{new}(\Gamma_0(8), \chi_8)$ and $S_4^{new}(\Gamma_0(24), \chi_8)$. A basis for $S_4^{new}(\Gamma_0(8),\chi_8)$ is given in \eqref{8chi8}.  The following six eta-quotients span the 
space $S_4^{new}(\Gamma_0(24),\chi_8)$.
\begin{equation}
\begin{split}\label{24:chi8}
f_{4,24,\chi_8;1}(z) &=  1^2 2^1 3^{-4} 4^1 6^{10} 8^2 12^{-4} :=  \sum_{n\ge 1} a_{4,24,\chi_8;1}(n) q^n ,~~\\
f_{4,24,\chi_8;2}(z) & =  1^1 2^3 3^{-1} 4^1 6^4 8^{-1} 24^1 := \sum_{n\ge 1} a_{4,24,\chi_8;2}(n) q^n,~~\\
f_{4,24,\chi_8;3}(z) & =  1^{-1} 2^4 3^1 6^3 8^1 12^1 24^{-1} := \sum_{n\ge 1} a_{4,24,\chi_8;3}(n) q^n , ~~\\
f_{4,24,\chi_8;4}(z) & = 1^{-2} 2^4 4^2 6^1 8^2 12^1  := \sum_{n\ge 1} a_{4,24,\chi_8;4}(n) q^n,\\
f_{4,24,\chi_8;5}(z) & = 2^1 3^{-2} 4^1 6^4 12^2 24^2 := \sum_{n\ge 1} a_{4,24,\chi_8;5}(n) q^n , ~~\\
f_{4,24,\chi_8;6}(z) & = 1^{-6} 2^{14} 6^1 8^{-2} 12^1  := \sum_{n\ge 1} a_{4,24,\chi_8;6}(n) q^n
\end{split}
\end{equation}
A basis for the space $M_4(\Gamma_0(48),\chi_{8})$ is given in the following proposition. 

\begin{prop}\label{chi8} 
A basis for the space $M_4(\Gamma_0(48),\chi_8)$ is given by 
\begin{equation}
\begin{split}
&\left\{E_{4, {\bf 1},\chi_8}(tz),~ E_{4,\chi_8, {\bf 1}}(tz), t\vert 6; f_{4,8,\chi_8;1}(t_1z),f_{4,8,\chi_8;2}(t_1z), t_1\vert 6;  f_{4,24,\chi_8;1}(t_2z),\right.\\
& \left.f_{4,24,\chi_8;2}(t_2z),f_{4,24,\chi_8;3}(t_2z),f_{4,24,\chi_8;4}(t_2z),  f_{4,24,\chi_8;5}(t_2z),f_{4,24,\chi_8;6}(t_2z), t_2\vert 2\right\},\\
\end{split}
\end{equation}
where $E_{4,{\bf 1},\chi_8}(z)$ and $E_{4,\chi_8,{\bf 1}}(z)$ are defined in \eqref{eis:8chi8}, $f_{4,8,\chi_8;i}(z)$, $i=1,2$ are defined in \eqref{8chi8} and 
$f_{4,24,\chi_8;j}(z)$, $1\le j\le 6$ are defined by \eqref{24:chi8}
\end{prop}
For the sake of simplifying the notation, we shall list the basis in \propref{chi8} as \\
$\left\{g_\alpha(z)\vert 1\le \alpha\le 28\right\}$, where 
$g_1(z) = E_{4, {\bf 1},\chi_8}(z)$, $g_2(z) = E_{4, {\bf 1},\chi_8}(2z)$,
$g_3(z) = E_{4, {\bf 1},\chi_8}(3z)$,  $g_4(z) = E_{4, {\bf 1},\chi_8}(6z)$, 
$g_5(z) = E_{4,\chi_8, {\bf 1}}(z)$,  $g_6(z) = E_{4,\chi_8, {\bf 1}}(2z)$,  
$g_7(z) = E_{4,\chi_8, {\bf 1}}(3z)$, $g_8(z) = E_{4,\chi_8, {\bf 1}}(6z)$, 
$g_9(z) = f_{4,8,\chi_8;1}(z)$, $g_{10}(z) = f_{4,8,\chi_8;1}(2z)$, $g_{11}(z) = 
f_{4,8,\chi_8;1}(3z)$, $g_{12}(z) = f_{4,8,\chi_8;1}(6z)$, $g_{13}(z) = f_{4,8,\chi_8;2}(z)$,  
$g_{14}(z) = f_{4,8,\chi_8;2}(2z)$, $g_{15}(z) = f_{4,8,\chi_8;2}(3z)$, $g_{16}(z) = 
f_{4,8,\chi_8;2}(6z)$, 
$g_{17}(z) =$ 
$f_{4,24,\chi_8;1}(z)$, $g_{18}(z) = f_{4,24,\chi_8;1}(2z)$,
$g_{19}(z) = f_{4,24,\chi_8;2}(z)$, 
$g_{20}(z) = f_{4,24,\chi_8;2}(2z)$, 
$g_{21}(z) = f_{4,24,\chi_8;3}(z)$, 
$g_{22}(z) = f_{4,24,\chi_8;3}(2z)$, 
$g_{23}(z) = f_{4,24,\chi_8;4}(z)$,  
$g_{24}(z) = f_{4,24,\chi_8;4}(2z)$,
$g_{25}(z) = f_{4,24,\chi_8;5}(z)$, 
\break 
$g_{26}(z) = f_{4,24,\chi_8;5}(2z)$,
$g_{27}(z) = f_{4,24,\chi_8;6}(z)$, 
$g_{28}(z) = f_{4,24,\chi_8;6}(2z)$. 

\noindent 
As before, we also write the Fourier expansions of these basis elements as $g_\alpha(z) = \displaystyle{\sum_{n\ge 1}} B_\alpha(n) q^n$, $1\le \alpha\le 28$.

In this case, all the 10 quadruples  corresponding to the $\chi_8$ character space (in Table 1) have the property that the powers of the theta functions corresponding to the coefficients $2$ are odd and corresponding to $3$ are even. Therefore, the resulting products of theta functions are modular forms of weight $4$ on $\Gamma_0(48)$ with character $\chi_8$ (we use Fact II to prove this). So, we can express these products of theta functions as a linear combination of the basis given in \propref{chi8} as follows. 
\begin{equation}
\Theta^i(z) \Theta^j(2z)\Theta^k(3z)\Theta^l(4z) = \sum_{\alpha=1}^{28} b_\alpha g_\alpha(z).
\end{equation}
Comparing the $n$-th Fourier coefficients on both the sides, we get 
\begin{equation*}
N(1^i,2^j,3^k,4^l;n) = \sum_{\alpha=1}^{28} b_\alpha B_\alpha(n).
\end{equation*}
Explicit values for the constants $b_\alpha$, $1\le \alpha\le 28$ corresponding to these 10 cases are given in Table 4.

\subsection{A basis for $M_4(\Gamma_0(48), \chi_{12})$ and proof of \thmref{2}(iii).} 

The dimension of the space in this case is $30$, with $\dim_{\mathbb C}{\mathcal E}_4(\Gamma_0(48), \chi_{12}) = 12$ and 
$\dim_{\mathbb C}S_4(\Gamma_0(48), \chi_{12}) = 18$. Regarding the old class, the space $S_4^{new}(\Gamma_0(12), \chi_{12})$ has dimension 
$4$ and is spanned by the following four eta-quotients:
\begin{equation}
\begin{split}
f_{4,12,\chi_{12};1}(z) ~=~ 2^{-1} 3^4 4^2 6^5 12^{-2}, & \quad  
f_{4,12,\chi_{12};2}(z) ~=~ 
3^4 4^3 6^{-2} 12^3, \\
f_{4,12,\chi_{12};3}(z) ~=~ 2^2 3^4 4^{-1} 6^{-4} 12^7, & \quad 
f_{4,12,\chi_{12};4}(z) = 
1^4 4^{-1} 6^{-2} 12^7. \\
\end{split}
\end{equation}
We write the Fourier expansions of these forms as $f_{4,12,\chi_{12};j}(z) = \displaystyle{\sum_{n\ge 1}} a_{4,12,\chi_{12};j}(n) q^n$, $1\le j\le 4$. 
In order to get the span of the newforms space of level $48$, we need the following three eta-quotients of level $48$ and character $\chi_{12}$.

\begin{equation}
\begin{split}
& 1^{-4} 2^7 4^5 6^{-3} 8^{-3} 12^9 24^{-3} := \sum_{n\ge 1} a_{4,48,\chi_{12};1}(n) q^n, ~~ 
2^{-3} 3^4 4^9 6^{7} 8^{-3} 12^5 24^{-3} := \sum_{n\ge 1} a_{4,48,\chi_{12};2}(n) q^n, \\
& 1^{-2} 2^2 3^2 4^2 8^{1} 12^2 24^{1} := \sum_{n\ge 1} a_{4,48,\chi_{12};3}(n) q^n.  
\end{split}
\end{equation}
Using the above three eta-quotients, we obtain the span of the space $S_4^{new}(\Gamma_0(48),\chi_{12})$ given by 
the following six forms.
\begin{equation}
\begin{split}
f_{4,48,\chi_{12};1}(z) =\sum_{n\ge 1 \atop{n\equiv 1 (\textrm{mod} 4)}} a_{4,48,\chi_{12};1}(n) q^n ,\quad 
f_{4,48,\chi_{12};2}(z) =\sum_{n\ge 1 \atop{n\equiv  3 (\textrm{mod} 4)}} a_{4,48,\chi_{12};1}(n) q^n, &\\
f_{4,48,\chi_{12};3}(z) =\sum_{n\ge 1 \atop{n\equiv  1 (\textrm{mod} 4)}} a_{4,48,\chi_{12};2}(n) q^n, \quad 
f_{4,48,\chi_{12};4}(z) =\sum_{n\ge 1 \atop{n\equiv  3 (\textrm{mod} 4)}} a_{4,48,\chi_{12};2}(n) q^n, & \\
f_{4,48,\chi_{12};5}(z) =\sum_{n\ge 1 \atop{n\equiv  1 (\textrm{mod} 4)}} a_{4,48,\chi_{12};3}(n) q^n, \quad 
f_{4,48,\chi_{12};6}(z) =\sum_{n\ge 1 \atop{n\equiv  3 (\textrm{mod} 4)}} a_{4,48,\chi_{12};3}(n) q^n, &\\
\end{split}
\end{equation}
A basis for the space $(M_4(\Gamma_0(48)),\chi_{12})$ is given in the following proposition. 

\begin{prop}\label{chi12}
A basis for the space $M_4(\Gamma_0(48),\chi_{12})$ is given by 
\begin{equation}
\begin{split}
\left\{E_{4, {\bf 1},\chi_{12}}(tz),{E_{4,\chi_{12}, {\bf 1}}(tz),E_{4,\chi_{-4},\chi_{-3}}(tz),E_{4,\chi_{-3},\chi_{-4}}(tz)}, t\vert 4; 
f_{4,12,\chi_{12};j}(t_1z), t_1\vert 4, 1\le j\le 4; \right.&\\
\left. f_{4,48,\chi_{12};1}(z),  f_{4,48,\chi_{12};2}(z),f_{4,48,\chi_{12};3}(z), f_{4,48,\chi_{12};4}(z), f_{4,48,\chi_{12};5}(z),f_{4,48,\chi_{12};6}(z)\right\},&\\
\end{split}
\end{equation}
where the Eisenstein series in the basis are defined by \eqref{eisenstein}.
\end{prop}
Let us denote the 30 basis elements in the above proposition as follows. \\
$\left\{h_\alpha(z)\vert 1\le \alpha\le 30\right\}$, where $h_1(z) = E_{4, {\bf 1},\chi_{12}}(z)$, $h_2(z) = E_{4,\chi_{12}, {\bf 1}}(z)$, 
$h_3(z) = E_{4,\chi_{-4},\chi_{-3}}(z)$, $h_4(z) = E_{4,\chi_{-3},\chi_{-4}}(z)$,  $h_5(z) = E_{4, {\bf 1},\chi_{12}}(2z)$, $h_6(z) = E_{4,\chi_{12}, {\bf 1}}(2z)$, 
$h_7(z) = E_{4,\chi_{-4},\chi_{-3}}(2z)$, $h_8(z) = 
\break 
E_{4,\chi_{-3},\chi_{-4}}(2z)$, $h_9(z) = E_{4, {\bf 1},\chi_{12}}(4z)$, $h_{10}(z) = 
E_{4,\chi_{12}, {\bf 1}}(4z)$, $h_{11}(z) = E_{4,\chi_{-4},\chi_{-3}}(4z)$, $h_{12}(z) = E_{4,\chi_{-3},\chi_{-4}}(4z)$, $h_{12+j}(z) =  f_{4,12,\chi_{12};j}(z)$, 
$1\le j\le 4$, $h_{16+j}(z) =  f_{4,12,\chi_{12};j}(2z)$, $1\le j\le 4$, $h_{20+j}(z) =  f_{4,12,\chi_{12};j}(4z)$, $1\le j\le 4$, $h_{24+j}(z) =  
f_{4,48,\chi_{12};j}(z)$, $1\le j\le 6$. 

\smallskip

To prove \thmref{2}(iii), we consider the case of  10 quadruples  corresponding to the $\chi_{12}$ character space (in Table 1). Now the roles of coefficients have interchanged and they have the property that the powers of the theta functions corresponding to the coefficients $2$ are even and corresponding to $3$ are odd. Therefore, the resulting products of theta functions are modular forms of weight $4$ on $\Gamma_0(48)$ with character $\chi_{12}$ (once again we  use Fact II to get this). So, we can express these products of theta functions as a linear combination of the basis given in \propref{chi12} as follows. 
\begin{equation}
\Theta^i(z) \Theta^j(2z)\Theta^k(3z)\Theta^l(4z) = \sum_{\alpha=1}^{30} c_\alpha h_\alpha(z),
\end{equation}
where $c_\alpha$'s are some constants. By comparing the $n$-th Fourier coefficients on both the sides, we get 
\begin{equation*}
N(1^i,2^j,3^k,4^l;n) = \sum_{\alpha=1}^{30} c_\alpha C_\alpha(n),
\end{equation*}
where  $h_\alpha(z) = \sum_{n\ge 1} C_\alpha(n) q^n$, $1\le \alpha\le 30$. Explicit values for the constants $c_\alpha$, $1\le \alpha\le 30$ corresponding to these 
10 cases are given in Table 5. 

\subsection{A basis for $M_4(\Gamma_0(48), \chi_{24})$ and proof of \thmref{2}(iv).} ~ In this case the space has dimension $28$. 
To get the span of the Eisenstein series space ${\mathcal E}_4(\Gamma_0(48),\chi_{24})$  (which has dimension 8), 
we use the Eisenstein series $E_{4,\chi,\psi}(z)$ defined in \eqref{eisenstein}, where 
$\chi,\psi \in \{{\bf 1}, \chi_{-8}, \chi_{-12}, \chi_{24}\}$. Note that $S_4^{new}(\Gamma_0(48),\chi_{24}) = \{0\}$.  
The space $S_4^{new}(\Gamma_0(24),\chi_{24})$ is spanned by the following ten eta-quotients (notation as in \eqref{eta-q}):
\begin{equation}
\begin{split}
& f_{4,24,\chi_{24};1}(z) ~=~ 3^{-2} 6^7 8^3 12^3 24^{-3}, ~
f_{4,24,\chi_{24};2}(z) ~=~ 3^{2} 4^7 6^{-3} 8^{-2} 12^{4}, ~
f_{4,24,\chi_{24};3}(z) = 3^{2} 4^{-3} 6^{1} 8^{6} 12^{2}, \\
& f_{4,24,\chi_{24};4}(z) =  3^2 6^{-3} 8^3 12^5 24^1, ~\quad 
f_{4,24,\chi_{24};5}(z) = 3^2 4^2 6^{-3} 8^{-1} 12^3 24^5, ~\\
& f_{4,24,\chi_{24};6}(z) = 3^2 4^1 6^1 8^{-2} 12^{-2} 24^8, \quad  
 f_{4,24,\chi_{24};7}(z) = 3^2 4^1 6^1 8^{-2} 12^{-2} 24^8, ~\\
& f_{4,24,\chi_{24};8}(z) = 1^1 3^{-1} 6^1 8^{-2} 12^1 24^8, 
 f_{4,24,\chi_{24};9}(z) = 2^2 3^6 4^1 6^{-3} 8^2, ~
f_{4,24,\chi_{24};10}(z) = 3^2 4^3 6^5 12^{-4} 24^2. \\
\end{split}
\end{equation}
We write the Fourier expansions as $f_{4,24,\chi_{24};j}(z) = \sum_{n \ge 1} a_{4,24,\chi_{24};j}(n) q^n$.  We now give a basis for the space 
$M_4(\Gamma_0(48),\chi_{24})$ in the following proposition.

\begin{prop}\label{chi24}
The following functions span the space $M_4(\Gamma_0(48),\chi_{24})$. 
\begin{equation}
\begin{split}
& \quad \left\{E_{4, {\bf 1},\chi_{24}}(tz),{E_{4,\chi_{24}, {\bf 1}}(tz),E_{4,\chi_{-8},\chi_{-3}}(tz),E_{4,\chi_{-3},\chi_{-8},}(tz)}, t\vert 2; \right. \\   
& \left. f_{4,24,\chi_{24};j}(z), 1\le j\le 10; f_{4,24,\chi_{24};j}(2z), 1\le j\le 10\right\}.  
\end{split}
\end{equation}
\end{prop}
We list these basis elements as $\{F_\alpha(z)\vert 1\le \alpha\le 28\}$, where  
$ F_1(z) = E_{4, {\bf 1},\chi_{24}}(z)$, $ F_2(z) = E_{4, {\bf 1},\chi_{24}}(2z)$, 
$ F_3(z) =E_{4,\chi_{-8},\chi_{-3}}(z)$,  $F_4(z) = E_{4,\chi_{-8},\chi_{-3}}(2z)$, 
$ F_5(z) = E_{4,\chi_{24}, {\bf 1}}(z)$, $ F_6(z) = E_{4,\chi_{24}, {\bf 1}}(2z) $, 
$ F_7(z) = E_{4,\chi_{-3},\chi_{-8},}(z)$, $ F_8(z) = E_{4,\chi_{-3},\chi_{-8},}(2z)$, 
$F_9(z) =f_{4,24,\chi_{24};1}(z)$, $ F_{10}(z) = f_{4,24,\chi_{24};1}(2z)$, $ F_{11}(z) = f_{4,24,\chi_{24};2}(z)$, 
$ F_{12}(z) = f_{4,24,\chi_{24};2}(2z)$, $ F_{13}(z) = f_{4,24,\chi_{24};3}(z)$, $ F_{14}(z) =  f_{4,24,\chi_{24};3}(2z)$, $ F_{15}(z) = f_{4,24,\chi_{24};4}(z)$, 
$ F_{16}(z) = f_{4,24,\chi_{24};4}(2z)$,  $ F_{17}(z) = f_{4,24,\chi_{24};5}(z)$, $ F_{18}(z) = f_{4,24,\chi_{24};5}(2z)$, $ F_{19}(z) = f_{4,24,\chi_{24};6}(z)$, 
$ F_{20}(z) = f_{4,24,\chi_{24};6}(2z)$, $ F_{21}(z) = f_{4,24,\chi_{24};7}(z)$, $ F_{22}(z) = f_{4,24,\chi_{24};7}(2z)$, $ F_{23}(z) =f_{4,24,\chi_{24};8}(z)$, $ F_{24}(z) =f_{4,24,\chi_{24};8}(2z)$, $ F_{25}(z) = f_{4,24,\chi_{24};9}(z)$, 
$ F_{26}(z) = f_{4,24,\chi_{24};9}(2z)$, $ F_{27}(z) = f_{4,24,\chi_{24};10}(z)$, $ F_{28}(z) = f_{4,24,\chi_{24};10}(2z)$.

\noindent 
We express the Fourier coefficients of the function $F_\alpha(z)$ as $\sum_{n\ge 1} D_\alpha(n) q^n$, $1\le \alpha\le 28$. 
In this case all the 10 quadruples have the property that the powers of the theta functions corresponding to the coefficients $2$ and $3$ are both odd.
Therefore, the resulting functions belong to the space of modular forms of weight $4$ on $\Gamma_0(48)$ with character $\chi_{24}$. As before, 
we use Fact II to arrive at this result.  So, one can express these theta functions as a linear combination of the basis elements:
\begin{equation}
\Theta^i(z) \Theta^j(2z)\Theta^k(3z)\Theta^l(4z) = \sum_{\alpha=1}^{28} d_\alpha F_\alpha(z),
\end{equation}
where $d_\alpha$'s are some constants. Comparing the $n$-th Fourier coefficients on both the sides, we get 
\begin{equation*}
N(1^i,2^j,3^k,4^l;n) = \sum_{\alpha=1}^{28} d_\alpha D_\alpha(n).
\end{equation*}
The values for the constants $d_\alpha$, $1\le \alpha\le 28$ corresponding to the 10 cases are given in Table 6.

\smallskip

\noindent {\bf Remark 4.1.} As we see in the proofs of Theorems 2.1 and 2.2, the main part of the proof is the construction of a basis of modular forms of weight $4$ on $\Gamma_0(48)$ with some character.  In the case of Theorem 2.3 we need to consider the quadratic forms given by Table 2. Since all these quadratic forms correspond to modular forms of weight $4$ on $\Gamma_0(48)$ with character (as specified in the table), we can follow the proof of Theorem 2.2 and complete the proof. The corresponding coefficients are given in Tables 7 to 10. So, we omit the details here. 

\smallskip

\noindent {\bf Remark 4.2.} To prove our theorems, we determined explicit bases for the 
spaces of modular forms $M_4(\Gamma_0(48),\chi)$, where $\chi$ is either trivial character 
modulo $48$ or $\chi=\chi_m$, $m=8, 12$ or $24$.  For our construction of bases we used the 
newforms theory which involves old classes of levels $6, 8, 12, 16$ and $24$. Therefore, from our construction it is easy to get bases for the modular form spaces $M_4(\Gamma_0(24),\chi)$, 
where $\chi$ is either the trivial character modulo $24$ or $\chi=\chi_m$, $m=8, 12$ or $24$. 
In our work \cite{alladi-conf}, we used the above bases of $M_4(\Gamma_0(24), \chi)$ to find 
formulas for the representation numbers of a positive integer by certain classes of quadratic forms in $8$ variables. 

\smallskip

\noindent{\bf Remark 4.3} We now make a remark about the number of coefficients of the products of the theta functions that are needed to get the 
linear combination constants $a_\alpha$, $b_\alpha$, $c_\alpha$, $d_\alpha$ (and $a_\alpha'$, $b_\alpha'$, $c_\alpha'$, $d_\alpha'$) that are 
listed in the tables (3 to 6) (resp. 7 to 10). The Sturm bound for the determination 
of a non-zero modular form of level $N$ is $k\nu/12$, where $\nu$ is the index of the congruence subgroup $\Gamma_0(N)$ in 
$SL_2({\mathbb Z})$ (see \cite{sturm}). In our case $N=48$, $k=4$, the index $\nu = 96$ and so the Sturm bound is $32$. 
In the table below, we list the number of coefficients needed (to get the constants), the Sturm bound and the dimension of the respective 
spaces of modular forms. \\


\begin{center}
\begin{tabular}{|c|c|c|c|c|}
\hline 
Space & no. of coefficients & Sturm bound & dimension& Theorem\\
&needed&&&\\
\hline 
$M_4(\Gamma_0(48))$ & 33 & 32 & 30 & 2.2 (i) and 2.3 (i)\\
\hline 
$M_4(\Gamma_0(48), \chi_8)$ & 28 & 32 & 28 & 2.2 (ii) and 2.3 (ii)\\
\hline 
$M_4(\Gamma_0(48), \chi_{12})$ & 31 & 32 & 30 & 2.2 (iii) and 2.3 (iii)\\
\hline 
$M_4(\Gamma_0(48), \chi_{24})$ & 28 & 32 & 28 & 2.2 (iv) and 2.3 (iv)\\
\hline 
\end{tabular}
\end{center}

\bigskip

\noindent {\bf Remark 4.4} In a recent result, Z. S. Aygin \cite{{aygin1}, {aygin2}} proved an orthogonal relation and used it to compute the coefficients of Eisenstein series part of a given modular form $f$ of weight $2k$ on $\Gamma_0(N)$ ($N$ is odd and square-free) in terms of sum of divisors function. In particular for the modular form given by the products of the theta function $T(z) = \prod_{\delta\vert N} 
\Theta(\delta z)^{r_\delta}$, his result yields explicit constants of the Eisenstein series part of the function $T(z)$ which involve the exponents $r_\delta$. For the interested reader, we refer to \cite[Theorem 7.1.2]{aygin1}, \cite[Theorem 5.2]{aygin2}. We also refer to \cite[Theorem 2.1]{cooper}. 

\subsection{Sample formulas}
 In this section we shall give explicit formulas for Theorems 2.2 and 2.3 for a few cases. 

\smallskip

\noindent {\bf  Formulas for the cases $(5,0,2,1), (3,2,2,1)$ of \thmref{2}(i).} 
\begin{equation*}
\begin{split}
&{N(1^5,3^2,4^1;n)}  = \hskip 13.5cm \\ 
&\frac{14}{5}   \sigma_{3}(n)-\frac{54}{5} \sigma_{3}(n/3)- \frac{238}{5} \sigma_3(n/4) + \frac{252}{5} \sigma_{3}(n/8) +\frac{918}{5} \sigma_{3}(n/12)-\frac{448}{5} \sigma_{3}(n/16)-\\
&\frac{972}{5} \sigma_{3}(n/24)+\frac{1728}{5} \sigma_{3}(n/48)+\frac{7}{20}   \sigma_{3;\chi_{-4},\chi_{-4}}(n)+ \frac{27}{20}   \sigma_{3;\chi_{-4},\chi_{-4}}(n/3)- \frac{4}{5} a_{4,6}(n)+ \frac{16}{5} a_{4,6}(n/2) -\\
&\frac{176}{5} a_{4,6}(n/4)-\frac{1408}{5} a_{4,6}(n/8) + a_{4,8}(n/2)+ 27 a_{4,8}(n/6)+4 a_{4,12}(n)+ 20 a_{4,12}(n/4)+ \frac{1}{4} a_{4,16}(n)-\\
&\frac{27}{4} a_{4,16}(n/3)+9a_{4,24}(n/2)\left(\frac{4}{n/2}\right)+ \frac{22}{5}a_{4,6}(n)\left(\frac{-4}{n}\right)+ \frac{5}{4}a_{4,12}(n)\left(\frac{-4}{n}\right)-\frac{9}{4}a_{4,24}(n)\left(\frac{-4}{n}\right),
\end{split}
\end{equation*}


\begin{equation*}
\begin{split}
&{ N(1^3,2^2,3^2,4^1;n) }=\hskip 13.5cm \\
&\frac{7}{5} \sigma_{3}(n)-\frac{7}{5} \sigma_{3}(n/2)-\frac{27}{5} \sigma_{3}(n/3)+\frac{27}{5} \sigma_{3}(n/6)+\frac{28}{5} \sigma_{3}(n/8)-\frac{448}{5} \sigma_{3}(n/16) - \\
&\frac{108}{5} \sigma_{3}(n/24) +\frac{1728}{5} \sigma_{3}(n/48) - \frac{2}{5} a_{4,6}(n)- 24 a_{4,6}(n/4)- \frac{448}{5} a_{4,6}(n/8) - \frac{1}{2} a_{4,8}(n) -  a_{4,8}(n/2)-\\
&\frac{27}{2} a_{4,8}(n/3) - 27 a_{4,8}(n/6)+ 2 a_{4,12}(n)+2 a_{4,12}(n/2)+20 a_{4,12}(n/4)+ \frac{1}{2} a_{4,16}(n)-\\
&\frac{27}{2} a_{4,16}(n/3) + \frac{3}{2}a_{4,24}(n)\left(\frac{4}{n}\right)+3 a_{4,24}(n/2)\left(\frac{4}{n/2}\right)+ 3 a_{4,6}(n)\left(\frac{-4}{n}\right)- \frac{3}{2} a_{4,24}(n)\left(\frac{-4}{n}\right).\\
\end{split}
\end{equation*}

\smallskip

\noindent {\bf  Formulas for the cases $ (4,1,2,1), (2,3,2,1)$ of \thmref{2}(ii).} 
\begin{equation*}
\begin{split}
&{N(1^4,2^1,3^2,4^1;n)} =\hskip 13.5 cm \\
&-\frac{26}{451} \sigma_{3;{\bf 1},\chi_{8}}(n/2)+ \frac{108}{451} \sigma_{3;{\bf 1},\chi_{8}}(n/6)+ \frac{832}{451} \sigma_{3;\chi_{8},{\bf 1}}(n)+ \frac{3456}{451} \sigma_{3;\chi_{8},{\bf 1}}(n/3)- \frac{12448}{451} a_{4,8,\chi_8;1}(n)+\\
&\frac{38416}{451} a_{4,8,\chi_8;1}(n/2)+ \frac{1296}{451} a_{4,8,\chi_8;1}(n/3)+\frac{98496}{451} a_{4,8,\chi_8;1}(n/6)+ \frac{14016}{451} a_{4,8,\chi_8;2}(n)-\frac{28032}{451} a_{4,8,\chi_8;2}(n/2)-\\
&\frac{1728}{451} a_{4,8,\chi_8;2}(n/3)-\frac{3456}{451} a_{4,8,\chi_8;2}(n/6)+ \frac{668}{41} a_{4,24,\chi_8;1}(n)- \frac{2368}{41} a_{4,24,\chi_8;1}(n/2)+\frac{3240}{41} a_{4,24,\chi_8;2}(n)-\\
& \frac{11208}{41} a_{4,24,\chi_8;2}(n/2)-\frac{216}{41} a_{4,24,\chi_8;3}(n)+ \frac{2856}{41} a_{4,24,\chi_8;3}(n/2)- \frac{6240}{41} a_{4,24,\chi_8;4}(n)+ \frac{20672}{41} a_{4,24,\chi_8;4}(n/2)-\\
&\frac{11232}{41} a_{4,24,\chi_8;5}(n)+ \frac{46656}{41} a_{4,24,\chi_8;5}(n/2)+ \frac{932}{41} a_{4,24,\chi_8;6}(n)- \frac{3568}{41} a_{4,24,\chi_8;6}(n/2),\\
\end{split}
\end{equation*}


\begin{equation*}
\begin{split}
&{ N(1^2,2^3,3^2,4^1;n)} =\hskip 13.5 cm \\
&-\frac{26}{451} \sigma_{3;{\bf 1},\chi_{8}}(n/2)+ \frac{108}{451} \sigma_{3;{\bf 1},\chi_{8}}(n/6)+ \frac{416}{451} \sigma_{3;\chi_{8},{\bf 1}}(n)+ \frac{1728}{451} \sigma_{3;\chi_{8},{\bf 1}}(n/3)- \frac{6224}{451} a_{4,8,\chi_8;1}(n)+\\
&\frac{16768}{451} a_{4,8,\chi_8;1}(n/2)+ \frac{648}{451} a_{4,8,\chi_8;1}(n/3)+\frac{98496}{451} a_{4,8,\chi_8;1}(n/6)+ \frac{7008}{451} a_{4,8,\chi_8;2}(n)+\frac{832}{451} a_{4,8,\chi_8;2}(n/2)-\\
&\frac{864}{451} a_{4,8,\chi_8;2}(n/3)-\frac{3456}{451} a_{4,8,\chi_8;2}(n/6)+ \frac{416}{41} a_{4,24,\chi_8;1}(n)- \frac{1384}{41} a_{4,24,\chi_8;1}(n/2)+\frac{1620}{41} a_{4,24,\chi_8;2}(n)-\\
& \frac{6288}{41} a_{4,24,\chi_8;2}(n/2)-\frac{108}{41} a_{4,24,\chi_8;3}(n)+ \frac{1872}{41} a_{4,24,\chi_8;3}(n/2)- \frac{2464}{41} a_{4,24,\chi_8;4}(n)+ \frac{10176}{41} a_{4,24,\chi_8;4}(n/2)-\\
&\frac{5616}{41} a_{4,24,\chi_8;5}(n)+ \frac{23040}{41} a_{4,24,\chi_8;5}(n/2)+ \frac{384}{41} a_{4,24,\chi_8;6}(n)- \frac{1928}{41} a_{4,24,\chi_8;6}(n/2).\\
\end{split}
\end{equation*}

\smallskip

\noindent {\bf  Formulas for the cases $ (5,0,1,2), (3,2,1,2)$ of \thmref{2}(iii).} 

\begin{equation*}
\begin{split}
&{ N(1^5,3^1,4^2;n)} =\hskip 13.5 cm \\
&\frac{54}{23} \sigma_{3;\chi_{12},{\bf 1},}(n)-\frac{2}{23} \sigma_{3;\chi_{-4},\chi_{-3}}(n)+ \frac{108}{23} \sigma_{3;\chi_{12},{\bf 1},}(n/2)+\frac{4}{23} \sigma_{3;\chi_{-4},\chi_{-3}}(n/2)+ \frac{1}{23} \sigma_{3;{\bf 1},\chi_{12}}(n/4)-\\
&\frac{1728}{23} \sigma_{3;\chi_{12},{\bf 1},}(n/4)+\frac{64}{23} \sigma_{3;\chi_{-4},\chi_{-3}}(n/4)-\frac{27}{23} \sigma_{3;\chi_{-3},\chi_{-4}}(n/4)+ \frac{63}{23} a_{4,12,\chi_{12};1}(n)+ \frac{12}{23} a_{4,12,\chi_{12};2}(n)+\\
&\frac{54}{23} a_{4,12,\chi_{12};3}(n)- \frac{174}{23} a_{4,12,\chi_{12};4}(n)+ \frac{348}{23} a_{4,12,\chi_{12};1}(n/2)- \frac{924}{23} a_{4,12,\chi_{12};2}(n/2)-\frac{12}{23} a_{4,12,\chi_{12};3}(n/2)-\\
&\frac{120}{23} a_{4,12,\chi_{12};4}(n/2)+ \frac{632}{23} a_{4,12,\chi_{12};1}(n/4)- \frac{2528}{23} a_{4,12,\chi_{12};2}(n/4)+ \frac{2528}{23} a_{4,12,\chi_{12};3}(n/4)+\\
&\frac{-5056}{23} a_{4,12,\chi_{12};4}(n/4)+ 9 a_{4,48,\chi_{12};1}(n)+3 a_{4,48,\chi_{12};2}(n)- 4 a_{4,48,\chi_{12};3}(n)-64 a_{4,48,\chi_{12};5}(n),\\
\end{split}
\end{equation*}


\begin{equation*}
\begin{split}
&{ N(1^3,2^2,3^1,4^2;n)} =\hskip 13.5 cm \\
&\frac{27}{23} \sigma_{3;\chi_{12},{\bf 1},}(n)-\frac{1}{23} \sigma_{3;\chi_{-4},\chi_{-3}}(n)+ \frac{1}{23} \sigma_{3;{\bf 1},\chi_{12}}(n/4)-\frac{27}{23} \sigma_{3;\chi_{-3},\chi_{-4}}(n/4)+ \frac{43}{23} a_{4,12,\chi_{12};1}(n)- \frac{40}{23} a_{4,12,\chi_{12};2}(n)+\\
&\frac{73}{23} a_{4,12,\chi_{12};3}(n)- \frac{179}{23} a_{4,12,\chi_{12};4}(n)+8a_{4,12,\chi_{12};1}(n/2)- 32 a_{4,12,\chi_{12};2}(n/2)+ 32 a_{4,12,\chi_{12};3}(n/2)+\\
&8 a_{4,12,\chi_{12};4}(n/2)+ \frac{256}{23} a_{4,12,\chi_{12};1}(n/4)- \frac{1408}{23} a_{4,12,\chi_{12};2}(n/4)+ \frac{1312}{23} a_{4,12,\chi_{12};3}(n/4)- \frac{2528}{23} a_{4,12,\chi_{12};4}(n/4)+\\
&\frac{9}{2} a_{4,48,\chi_{12};1}(n)+ \frac{7}{2} a_{4,48,\chi_{12};2}(n)- \frac{3}{2} a_{4,48,\chi_{12};3}(n)+ \frac{3}{2} a_{4,48,\chi_{12};4}(n)- 24 a_{4,48,\chi_{12};5}(n)- 24 a_{4,48,\chi_{12};6}(n).\\
\end{split}
\end{equation*}


\noindent {\bf  Formulas for the cases $ (4,1,1,2), (2,3,1,2)$ of \thmref{2}(iv).} 
\begin{equation*}
\begin{split}
&{ N(1^4,2^1,3^1,4^2;n)}=\hskip 13.5 cm \\
&\frac{1}{261} \sigma_{3;{\bf 1},\chi_{24}}(n/2)+\frac{16}{261} \sigma_{3;\chi_{-8},\chi_{-3}}(n)+ \frac{48}{29} \sigma_{3;\chi_{24},{\bf 1}}(n)-\frac{3}{29} \sigma_{3;\chi_{-3},\chi_{-8}}(n/2)+ \frac{7432}{261} a_{4,24,\chi_{24};1}(n)+\\
&\frac{45424}{261} a_{4,24,\chi_{24};1}(n/2)+ \frac{1378}{87} a_{4,24,\chi_{24};2}(n)+ \frac{13576}{87} a_{4,24,\chi_{24};2}(n/2)- \frac{10214}{261} a_{4,24,\chi_{24};3}(n)-\\
&\frac{80672}{261} a_{4,24,\chi_{24};3}(n/2)- \frac{17048}{87} a_{4,24,\chi_{24};4}(n)- \frac{123680}{87} a_{4,24,\chi_{24};4}(n/2)- \frac{5804}{87} a_{4,24,\chi_{24};5}(n)-\\
&\frac{42944}{87} a_{4,24,\chi_{24};5}(n/2)+ \frac{57448}{261} a_{4,24,\chi_{24};6}(n)+ \frac{624352}{261} a_{4,24,\chi_{24};6}(n/2)+ \frac{14218}{261} a_{4,24,\chi_{24};7}(n)+\\
&\frac{252928}{261} a_{4,24,\chi_{24};7}(n/2)+ \frac{25472}{261} a_{4,24,\chi_{24};8}(n)- \frac{252928}{261} a_{4,24,\chi_{24};8}(n/2)- \frac{5792}{261} a_{4,24,\chi_{24};9}(n)-\\
&\frac{39656}{261} a_{4,24,\chi_{24};9}(n/2)- \frac{6418}{261} a_{4,24,\chi_{24};10}(n) - \frac{43888}{261} a_{4,24,\chi_{24};10}(n/2),\\
\end{split}
\end{equation*}

\begin{equation*}
\begin{split}
&{ N(1^2,2^3,3^1,4^2;n)}=\hskip 13.5 cm \\
&\frac{1}{261} \sigma_{3;{\bf 1},\chi_{24}}(n/2)+\frac{8}{261} \sigma_{3;\chi_{-8},\chi_{-3}}(n)+ \frac{24}{29} \sigma_{3;\chi_{24},{\bf 1}}(n)-\frac{3}{29} \sigma_{3;\chi_{-3},\chi_{-8}}(n/2)+\frac{5804}{261} a_{4,24,\chi_{24};1}(n)+\\
&  \frac{24544}{261} a_{4,24,\chi_{24};1}(n/2)+ \frac{1472}{87} a_{4,24,\chi_{24};2}(n)+\frac{9400}{87} a_{4,24,\chi_{24};2}(n/2)-\frac{9022}{261} a_{4,24,\chi_{24};3}(n)- \\
&\frac{51440}{261} a_{4,24,\chi_{24};3}(n/2)- \frac{15136}{87} a_{4,24,\chi_{24};4}(n)- \frac{72176}{87} a_{4,24,\chi_{24};4}(n/2)-\frac{4468}{87} a_{4,24,\chi_{24};5}(n)- \\
&\frac{13712}{87} a_{4,24,\chi_{24};5}(n/2)+ \frac{61088}{261} a_{4,24,\chi_{24};6}(n)+ \frac{377968}{261} a_{4,24,\chi_{24};6}(n/2)+\frac{19898}{261} a_{4,24,\chi_{24};7}(n)+ \\
&\frac{119296}{261} a_{4,24,\chi_{24};7}(n/2)+ \frac{21088}{261} a_{4,24,\chi_{24};8}(n)- \frac{119296}{261} a_{4,24,\chi_{24};8}(n/2)-\frac{22952}{261} a_{4,24,\chi_{24};9}(n/2)- \\
&\frac{5036}{261} a_{4,24,\chi_{24};10}(n)- \frac{25096}{261} a_{4,24,\chi_{24};10}(n/2).\\
\end{split}
\end{equation*}

\noindent {\bf  Formulas for the cases $ (0,0,4,2,2), (1,1,1,2,3)$ of \thmref{3}(i).} 
\begin{equation*}
\begin{split}
&{ N(3^4,4^2,6^2;n)} = \hskip 13.5cm \\
&  \frac{1}{10}\sigma_{3}(n)-\frac{1}{10} \sigma_{3}(n/2)-\frac{21}{10} \sigma_{3}(n/3)+\frac{21}{10} \sigma_{3}(n/6)+\frac{8}{10} \sigma_{3}(n/8) - \frac{64}{5} \sigma_{3}(n/16) - \frac{84}{5} \sigma_{3}(n/24)+\\ 
&\frac{1344}{25} \sigma_{3}(n/48)-\frac{4}{15} a_{4,6}(n)-\frac{4}{3} a_{4,6}(n/2)- 16 a_{4,6}(n/4)- \frac{384}{5} a_{4,6}(n/8)-\frac{1}{2} a_{4,8}(n) - 
5 a_{4,8}(n/2)-\\
&\frac{3}{2} a_{4,8}(n/3)-15 a_{4,8}(n/6)+\frac{1}{6} a_{4,12}(n)+2 a_{4,12}(n/2)+12 a_{4,12}(n/4) - a_{4,16}(n)+3 a_{4,16}(n/3)+ \\
&\frac{1}{2} a_{4,24}(n)\left(\frac{4}{n}\right)+ 3a_{4,24}(n/2)\left(\frac{4}{n/2}\right)+ \frac{4}{3}a_{4,6}(n)\left(\frac{-4}{n}\right)+\frac{2}{3} a_{4,12}(n)\left(\frac{-4}{n}\right)- a_{4,24}(n)\left(\frac{-4}{n}\right)\\
\end{split}
\end{equation*}

\begin{equation*}
\begin{split}
&{ N(1^1,2^1,3^1,4^2,6^3;n)} = \hskip 13.5cm \\
&\frac{1}{8}\sigma_{3}(n)-\frac{1}{8} \sigma_{3}(n/2)-\frac{9}{8} \sigma_{3}(n/3)+\frac{9}{8} \sigma_{3}(n/6)+2 \sigma_{3}(n/8)-32 \sigma_{3}(n/16) - 
18 \sigma_{3}(n/24)+288 \sigma_{3}(n/48)- \\
& 4 a_{4,6}(n/4)- 32 a_{4,6}(n/8)-\frac{1}{8} a_{4,8}(n)- \frac{5}{2}a_{4,8}(n/2)+\frac{9}{8} a_{4,8}(n/3) -
\frac{27}{2}a_{4,8}(n/6)+\frac{3}{8} a_{4,12}(n) + \\
&2 a_{4,12}(n/2)+6 a_{4,12}(n/4)-\frac{1}{4} a_{4,16}(n)-\frac{9}{4} a_{4,16}(n/3)+\frac{5}{8} a_{4,24}(n)\left(\frac{4}{n}\right) +
\frac{3}{2}a_{4,24}(n/2)\left(\frac{4}{n/2}\right) + \\
&a_{4,6}(n)\left(\frac{-4}{n}\right)+\frac{1}{2} a_{4,12}(n)\left(\frac{-4}{n}\right)- \frac{1}{4}a_{4,24}(n)\left(\frac{-4}{n}\right)\\
\end{split}
\end{equation*}

\noindent {\bf  Formulas for the cases $ (0,0,3,2,3), (0,0,1,4,3)$ of \thmref{3}(ii).} 
\begin{equation*}
\begin{split}
&{N(3^3,4^2,6^3;n)} =\hskip 13.5 cm \\ 
&\frac{4}{451} \sigma_{3;{\bf 1},\chi_{8}}(n/2)+ \frac{78}{451} \sigma_{3;{\bf 1},\chi_{8}}(n/6)+ \frac{32}{451} \sigma_{3;\chi_{8},{\bf 1}}(n)-\frac{624}{451} \sigma_{3;\chi_{8},{\bf 1}}(n/3)+ \frac{44746}{4059} a_{4,8,\chi_8;1}(n)-\\
&\frac{17944}{4059} a_{4,8,\chi_8;1}(n/2)- \frac{234}{451} a_{4,8,\chi_8;1}(n/3)-\frac{15456}{451} a_{4,8,\chi_8;1}(n/6)-\frac{7232}{451} a_{4,8,\chi_8;2}(n)-\frac{128}{451} a_{4,8,\chi_8;2}(n/2)+\\
&\frac{312}{451} a_{4,8,\chi_8;2}(n/3)-\frac{2496}{451} a_{4,8,\chi_8;2}(n/6)-\frac{2582}{369} a_{4,24,\chi_8;1}(n)+\frac{1664}{369} a_{4,24,\chi_8;1}(n/2)-\frac{3682}{123} a_{4,24,\chi_8;2}(n)+\\
& \frac{2776}{123} a_{4,24,\chi_8;2}(n/2)-\frac{650}{123} a_{4,24,\chi_8;3}(n)- \frac{1520}{123} a_{4,24,\chi_8;3}(n/2)+ \frac{23752}{369} a_{4,24,\chi_8;4}(n)- \frac{13888}{369} a_{4,24,\chi_8;4}(n/2)+\\
&\frac{4488}{41} a_{4,24,\chi_8;5}(n)-\frac{6976}{41} a_{4,24,\chi_8;5}(n/2)-\frac{1154}{41} a_{4,24,\chi_8;6}(n)+\frac{1016}{41} a_{4,24,\chi_8;6}(n/2),\\
\end{split}
\end{equation*}

\begin{equation*}
\begin{split}
&{N(3^1,4^4,6^3;n) }=\hskip 13.5 cm \\
&-\frac{8}{451} \sigma_{3;{\bf 1},\chi_{8}}(n/2)+ \frac{90}{451} \sigma_{3;{\bf 1},\chi_{8}}(n/6)+ \frac{16}{451} \sigma_{3;\chi_{8},{\bf 1}}(n)+ \frac{180}{451} \sigma_{3;\chi_{8},{\bf 1}}(n/3)- \frac{8177}{1353} a_{4,8,\chi_8;1}(n)+\\
&\frac{3544}{1353} a_{4,8,\chi_8;1}(n/2)-\frac{1962}{451} a_{4,8,\chi_8;1}(n/3)+\frac{17136}{451} a_{4,8,\chi_8;1}(n/6)+ \frac{3600}{451} a_{4,8,\chi_8;2}(n)+\frac{256}{451} a_{4,8,\chi_8;2}(n/2)+\\
&\frac{8028}{451} a_{4,8,\chi_8;2}(n/3)-\frac{2880}{451} a_{4,8,\chi_8;2}(n/6)+ \frac{499}{123} a_{4,24,\chi_8;1}(n)- \frac{344}{123} a_{4,24,\chi_8;1}(n/2)+\frac{589}{41} a_{4,24,\chi_8;2}(n)-\\
& \frac{648}{41} a_{4,24,\chi_8;2}(n/2)-\frac{83}{41} a_{4,24,\chi_8;3}(n)+ \frac{248}{41} a_{4,24,\chi_8;3}(n/2)- \frac{2984}{123} a_{4,24,\chi_8;4}(n)+ \frac{3136}{123} a_{4,24,\chi_8;4}(n/2)-\\
&\frac{2184}{41} a_{4,24,\chi_8;5}(n)+ \frac{3456}{41} a_{4,24,\chi_8;5}(n/2)+ \frac{163}{41} a_{4,24,\chi_8;6}(n)- \frac{240}{41} a_{4,24,\chi_8;6}(n/2).\\
\end{split}
\end{equation*}

\noindent {\bf  Formulas for the cases $ (0,0,5,1,2), (0,0,1,3,4)$ of \thmref{3}(iii).} 
\begin{equation*}
\begin{split}
&{N(3^5,4^1,6^2;n)}=\hskip 13.5 cm\\
&\frac{2}{23} \sigma_{3;\chi_{12},{\bf 1},}(n)+\frac{2}{23} \sigma_{3;\chi_{-4},\chi_{-3}}(n)+ \frac{1}{23} \sigma_{3;{\bf 1},\chi_{12}}(n/4)+\frac{1}{23} \sigma_{3;\chi_{-3},\chi_{-4}}(n/4)+ \frac{19}{23} a_{4,12,\chi_{12};1}(n)+\\
&\frac{22}{23} a_{4,12,\chi_{12};3}(n)+ \frac{22}{23} a_{4,12,\chi_{12};4}(n)+ 16 a_{4,12,\chi_{12};2}(n/2)+ 24 a_{4,12,\chi_{12};3}(n/2)- \frac{416}{23} a_{4,12,\chi_{12};1}(n/4)+\\
&\frac{2912}{23} a_{4,12,\chi_{12};3}(n/4)+ \frac{2912}{23} a_{4,12,\chi_{12};4}(n/4)- \frac{5}{2} a_{4,48,\chi_{12};1}(n)+ \frac{7}{6} a_{4,48,\chi_{12};2}(n)+\\
&\frac{3}{2} a_{4,48,\chi_{12};3}(n)+ \frac{31}{6} a_{4,48,\chi_{12};4}(n)+ \frac{40}{3} a_{4,48,\chi_{12};5}(n)- \frac{56}{3} a_{4,48,\chi_{12};6}(n),\\
\end{split}
\end{equation*}

\begin{equation*}
\begin{split}
&{N(3^1,4^3,6^4;n)}=\hskip 13.5 cm \\
&\frac{3}{92} \sigma_{3;\chi_{12},{\bf 1},}(n)-\frac{1}{92} \sigma_{3;\chi_{-4},\chi_{-3}}(n)+ \frac{1}{23} \sigma_{3;{\bf 1},\chi_{12}}(n/4)-\frac{3}{23} \sigma_{3;\chi_{-3},\chi_{-4}}(n/4)-\frac{1}{46} a_{4,12,\chi_{12};1}(n)-\\
& \frac{8}{23} a_{4,12,\chi_{12};2}(n) + \frac{173}{46} a_{4,12,\chi_{12};3}(n)+ \frac{17}{46} a_{4,12,\chi_{12};4}(n)+ 12 a_{4,12,\chi_{12};4}(n/2)+ \frac{140}{23} a_{4,12,\chi_{12};1}(n/4)+ \\
&\frac{272}{23} a_{4,12,\chi_{12};2}(n/4) -
\frac{1312}{23} a_{4,12,\chi_{12};3}(n/4)- \frac{1072}{23} a_{4,12,\chi_{12};4}(n/4)+ \frac{3}{4} a_{4,48,\chi_{12};1}(n)- \frac{1}{4} a_{4,48,\chi_{12};2}(n) -\\ 
&\frac{3}{4} a_{4,48,\chi_{12};3}(n) -
\frac{7}{4} a_{4,48,\chi_{12};4}(n)-4 a_{4,48,\chi_{12};5}(n)+ 4 a_{4,48,\chi_{12};6}(n).\\
\end{split}
\end{equation*}

\noindent {\bf  Formulas for the cases $ (0,0,4,1,3), (1,1,1,3,2)$ of \thmref{3}(iv).} 
\begin{equation*}
\begin{split}
&{N(3^4,4^1,6^3;n)}= \hskip 13.5cm \\
&\frac{1}{261}\sigma_{3;{\bf 1},\chi_{24}}(n/2)+\frac{16}{261} \sigma_{3;\chi_{-8},\chi_{-3}}(n)+ \frac{16}{261} \sigma_{3;\chi_{24},{\bf 1}}(n)-\frac{1}{261} \sigma_{3;\chi_{-3},\chi_{-8}}(n/2)+\\
&\frac{1504}{261} a_{4,24,\chi_{24};1}(n)+\frac{512}{87}a_{4,24,\chi_{24};1}(n/2)+\frac{512}{87}a_{4,24,\chi_{24};2}(n)+\frac{1320}{29} 
a_{4,24,\chi_{24};2}(n/2)-  \\
&\frac{1432}{261} a_{4,24,\chi_{24};3}(n) -\frac{8432}{87} a_{4,24,\chi_{24};3}(n/2)-\frac{10096}{261} a_{4,24,\chi_{24};4}(n)-\frac{132608}{261} a_{4,24,\chi_{24};4}(n/2)- \\
&\frac{1024}{87} a_{4,24,\chi_{24};5}(n)- 256 a_{4,24,\chi_{24};5}(n/2) +
\frac{6128}{87}a_{4,24,\chi_{24};6}(n)+\frac{227648}{261} a_{4,24,\chi_{24};6}(n/2)+\\
&\frac{1640}{87} a_{4,24,\chi_{24};7}(n)+\frac{3760}{9}a_{4,24,\chi_{24};7}(n/2) -
\frac{132608}{261}a_{4,24,\chi_{24};8}(n/2)-\frac{512}{87} a_{4,24,\chi_{24};9}(n)-\\
&\frac{16756}{261} a_{4,24,\chi_{24};9}(n/2)+\frac{512}{87} a_{4,24,\chi_{24};10}(n)-\frac{15536}{261} a_{4,24,\chi_{24};10}(n/2)\\
\end{split}
\end{equation*}
\begin{equation*}
\begin{split}
&{ N(1^1,2^1,3^1,4^3,6^2;n)} = \hskip 13.5cm \\
&\frac{1}{261}\sigma_{3;{\bf 1},\chi_{24}}(n/2)-\frac{4}{261} \sigma_{3;\chi_{-8},\chi_{-3}}(n)+ \frac{4}{29} \sigma_{3;\chi_{24},{\bf 1}}(n)+\frac{1}{29} \sigma_{3;\chi_{-3},\chi_{-8}}(n/2)+\\
&\frac{194}{261} a_{4,24,\chi_{24};1}(n)-\frac{13472}{261} a_{4,24,\chi_{24};1}(n/2)-\frac{202}{87} a_{4,24,\chi_{24};2}(n)-\frac{4664}{87} a_{4,24,\chi_{24};2}(n/2)+\\
&\frac{1190}{261} a_{4,24,\chi_{24};3}(n) + \frac{776}{9} a_{4,24,\chi_{24};3}(n/2)+\frac{1172}{87} a_{4,24,\chi_{24};4}(n)+\frac{1352}{3} 
a_{4,24,\chi_{24};4}(n/2)-\\
& \frac{304}{87} a_{4,24,\chi_{24};5}(n)+\frac{13144}{87}a_{4,24,\chi_{24};5}(n/2) +
\frac{164}{261} a_{4,24,\chi_{24};6}(n)-\frac{218696}{261} a_{4,24,\chi_{24};6}(n/2)+\\
&\frac{698}{261} a_{4,24,\chi_{24};7}(n)- \frac{68048}{261}a_{4,24,\chi_{24};7}(n/2)- \frac{14576}{261} a_{4,24,\chi_{24};8}(n)+
\frac{96128}{261}a_{4,24,\chi_{24};8}(n/2)+\\
&\frac{296}{261} a_{4,24,\chi_{24};9}(n)+\frac{13984}{261} a_{4,24,\chi_{24};9}(n/2)+ 
\frac{286}{261} a_{4,24,\chi_{24};10}(n)+\frac{12956}{261} a_{4,24,\chi_{24};10}(n/2)\\
\end{split}
\end{equation*}

\bigskip


\noindent {\bf Remark 4.5	.} We note that formulas for the cases $(i,j,k,l,m)$ (i.e., with  coefficients $1,2,3,4,6$) in Table 2, with $i=k=0$ can be obtained from the work \cite{a-k1} (replacing $n$ by $n/2$ in their formulas). There are 28 such cases (10 for $\chi_{24}$ character case and 6 each in the remaining 3 characters) in Table 2 with this condition. It is to be noted that different bases were used in \cite{a-k1} to get the formulas. So, replacing $n/2$ by $n$ in our formulas in these 28 cases, we get  different formulas for these cases (with coefficients 1,2,3) when compared with \cite{a-k1}.



\section{List of tables}
In this section we list Tables 1 and 2, which  give the list of exponents of the theta functions corresponding to the coefficients $1,2,3,4$ and 
$1,2,3,4,6$ respectively. We mention only the character $\chi$ in these tables corresponding to the space $M_4(\Gamma_0(48), \chi)$.  
The rest of the tables 3 to 10 give the coefficients list for the formulas for the number of representations corresponding to Theorem 2.2 and 
Theorem 2.3. These tables (3 to 10) are kept at the end as an Appendix.


\begin{center}
\textbf{Table 1.\\}
Octonary quadratic forms with coefficients $1,2,3,4$ ($k,l\not= 0$).
{\small 
\begin{tabular}{|c|c|}
\hline
$(i,j,k,l)$ & Space\\
\hline 
(0,0,6,2),(0,0,4,4),(0,0,2,6),(0,4,2,2),(0,2,4,2),(0,2,2,4),(5,0,2,1), &\\
(4,0,2,2),(3,0,4,1), (3,0,2,3),(2,0,4,2),(2,0,2,4),(1,0,6,1),(1,0,4,3), &\rm triv. \\
(1,0,2,5),(3,2,2,1),(2,2,2,2),(1,4,2,1),(1,2,4,1),(1,2,2,3)&\\
\hline 
(0,5,2,1),(0,3,4,1),(0,3,2,3),(0,1,6,1),(0,1,4,3),(0,1,2,5),(4,1,2,1), &\\
(3,1,2,2),(2,3,2,1),(2,1,4,1),(2,1,2,3),(1,3,2,2),(1,1,4,2),(1,1,2,4) &$\chi_{8}$\\
\hline
(0,0,7,1),(0,0,5,3),(0,0,3,5),(0,0,1,7),(0,6,1,1),(0,4,3,1),(0,4,1,3),(0,2,5,1), &\\
(0,2,3,3),(0,2,1,5),(6,0,1,1),(5,0,1,2),(4,0,3,1),(4,0,1,3),(3,0,3,2),(3,0,1,4), & $\chi_{12}$\\
(2,0,5,1),(2,0,3,3), (2,0,1,5),(1,0,5,2),(1,0,3,4),(1,0,1,6),(4,2,1,1),(3,2,1,2), &\\
(2,4,1,1),(2,2,3,1),(2,2,1,3),(1,4,1,2),(1,2,3,2), (1,2,1,4) &\\
\hline
(0,5,1,2),(0,3,3,2),(0,3,1,4),(0,1,5,2),(0,1,3,4),(0,1,1,6),(5,1,1,1), &\\
(4,1,1,2),(3,3,1,1).(3,1,3,1),(3,1,1,3),(2,3,1,2),(2,1,3,2),(2,1,1,4), & $\chi_{24}$\\
(1,5,1,1),(1,3,3,1),(1,3,1,3),(1,1,5,1),(1,1,3,3),(1,1,1,5)
&\\

\hline 
\end{tabular}
}
\end{center}

\begin{center}


\textbf{Table 2.}\\

Octonary quadratic forms with coefficients $1,2,3,4,6$ ($l,m\not= 0$).


{\small 
\begin{tabular}{|c|c|}
\hline
$(i,j,k,l,m)$ & Space\\
\hline 
(0,0,4,2,2),(0,0,2,4,2),(0,0,2,2,4),(0,2,2,2,2),(0,4,0,2,2), (0,2,0,4,2),(0,2,0,2,4),(0,1,5,1,1),&\\
(0,3,3,1,1),(0,5,1,1,1),(0,1,3,1,3),(0,3,1,1,3),(0,1,1,3,3), (0,1,1,1,5) (0,1,3,3,1),(0,3,1,3,1), &\\
(0,1,1,5,1),(1,1,1,2,3),(1,1,1,4,1),(1,1,3,2,1),(1,3,1,2,1),(1,0,4,1,2),(1,2,2,1,2),(1,4,0,1,2), &\\
(1,0,0,3,4),(1,0,2,3,2), (1,2,0,3,2),(1,0,0,5,2),(1,0,2,1,4), (1,2,0,1,4),(1,0,0,1,6),(2,2,0,2,2), &{\rm triv.}\\
(2,1,3,1,1),(2,3,1,1,1),(2,1,1,1,3),(2,1,1,3,1),(2,0,2,2,2), (2,0,0,4,2),(2,0,0,2,4),(3,1,1,2,1), &\\
(3,0,2,1,2),(3,2,0,1,2),(3,0,0,3,2),(3,0,0,1,4),(4,1,1,1,1),(4,0,0,2,2),(5,0,0,1,2)&
\\
\hline
(0,0,3,2,3),(0,0,1,4,3),(0,0,3,4,1),(0,0,5,2,1),(0,0,1,2,5),(0,0,1,6,1),(0,1,4,1,2),(0,3,2,1,2)&\\
(0,5,0,1,2),(0,1,0,3,4),(0,1,2,3,2),(0,3,0,3,2),(0,1,0,5,2),(0,1,2,1,4),(0,3,0,1,4),(0,1,0,1,6)&\\
(0,2,1,2,3),(0,2,1,4,1),(0,2,3,2,1),(0,4,1,2,1),(1,1,2,2,2),(1,1,0,4,2),(1,1,0,2,4),(1,3,0,2,2)&\\
(1,0,3,1,3),(1,2,1,1,3),(1,0,1,1,5),(1,0,1,3,3),(1,0,5,1,1),(1,2,3,1,1),(1,4,1,1,1),(1,0,3,3,1)&$\chi_8$\\
(1,2,1,3,1),(1,0,1,5,1),(2,2,1,2,1),(2,1,2,1,2),(2,3,0,1,2),(2,1,0,3,2),(2,1,0,1,4),(2,0,1,2,3)&\\
(2,0,1,4,1),(2,0,3,2,1),(3,1,0,2,2),(3,0,1,1,3),(3,0,3,1,1),(3,2,1,1,1),(3,0,1,3,1),(4,1,0,1,2)&\\
(4,0,1,2,1),(5,0,1,1,1)&\\
\hline
(0,0,5,1,2),(0,0,1,3,4),(0,0,3,3,2),(0,0,1,5,2),(0,0,3,1,4),(0,0,1,1,6),(0,2,3,1,2),(0,2,1,3,2)&\\
(0,4,1,1,2),(0,2,1,1,4),(0,3,0,2,3),(0,1,2,2,3),(0,1,0,4,3),(0,1,2,4,1),(0,3,0,4,1),(0,1,4,2,1)&\\
(0,3,2,2,1),(0,5,0,2,1),(0,1,0,2,5),(0,1,0,6,1),(1,1,2,1,3),(1,1,0,3,3),(1,1,4,1,1),(1,1,0,1,5)&\\
(1,1,2,3,1),(1,1,0,5,1),(1,0,3,2,2),(1,2,1,2,2)(1,0,1,4,2),(1,0,1,2,4),(1,3,0,1,3),(1,3,2,1,1)&$\chi_{12}$\\
(1,5,0,1,1),(1,3,0,3,1),(2,2,1,1,2),(2,0,3,1,2),(2,0,1,3,2),(2,0,1,1,4),(2,1,0,2,3),(2,1,0,4,1)&\\
(2,1,2,2,1),(2,3,0,2,1),(3,3,0,1,1),(3,0,1,2,2),(3,1,0,1,3),(3,1,2,1,1),(3,1,0,3,1),(4,0,1,1,2)&\\
(4,1,0,2,1),(5,1,0,1,1)&\\ 
\hline 
(0,0,4,1,3),(0,0,2,3,3),(0,0,6,1,1),(0,0,2,1,5), (0,0,4,3,1),(0,0,2,5,1),(0,2,2,1,3),(0,4,0,1,3), &\\
(0,2,0,3,3),(0,2,4,1,1),(0,4,2,1,1),(0,6,0,1,1), (0,2,0,1,5),(0,2,2,3,1),(0,4,0,3,1),(0,2,0,5,1), &\\
(0,1,3,2,2),(0,3,1,2,2),(0,1,1,4,2),(0,1,1,2,4), (1,1,1,3,2),(1,1,1,1,4),(1,1,3,1,2),(1,0,2,2,3), &\\
(1,2,0,2,3),(1,0,0,4,3),(1,0,2,4,1),(1,2,0,4,1), (1,0,4,2,1),(1,2,2,2,1),(1,4,0,2,1),(1,0,0,2,5), &$\chi_{24}$\\
(1,0,0,6,1),(1,3,1,1,2),(2,2,2,1,1),(2,2,0,1,3), (2,2,0,3,1),(2,0,2,1,3),(2,0,0,3,3),(2,0,4,1,1), &\\
(2,4,0,1,1),(2,0,0,1,5),(2,0,2,3,1),(2,0,0,5,1),(2,1,1,2,2),(3,0,0,2,3),(3,0,0,4,1),(3,0,2,2,1), &\\
(3,2,0,2,1),(3,1,1,1,2),(4,0,0,1,3),(4,0,2,1,1),(4,2,0,1,1),(4,0,0,3,1),(5,0,0,2,1),(6,0,0,1,1)&\\
\hline 
\end{tabular}
}
\end{center}


\section*{Acknowledgements}
We have used the open-source mathematics software SAGE (www.sagemath.org) to perform our calculations. Part of the work was done when the last named author visited the School of Mathematical Sciences, NISER, Bhubaneswar. He thanks the department for their hospitality and support. 
The second author is partially funded by SERB grant SR/FTP/MS-053/2012. We thank the referee for suggesting to include the data on comparing with the Sturm bound for finding the constants that appear in Tables 3 to 10 (Remark 4.3) and also for pointing out the works of Z. S. Aygin \cite{{aygin1}, {aygin2}} which give explicit values of the constants corresponding to the Eisenstein series part while expressing the products of the theta functions (Remark 4.4).

\bigskip

\section*{Appendix}

We list the tables 3 to 10 below.


\vskip 1cm 

\begin{center}
{\tiny
\textbf{Table 3.} (Theorem 2.2 (i))
\begin{tabular}{|c|c|c|c|c|c|c|c|c|c|c|c|c|c|c|c|}
\hline
$ijkl$&$a_1$ &$a_2$& $a_3$& $a_4$&$a_5$ &$a_6$& $a_7$& $a_8$&$a_9$ &$a_{10}$& $a_{11}$& $a_{12}$&$a_{13}$ &$a_{14}$& $a_{15}$ \\
\hline
&&&&&&&&&&&&&&&\\
$0062$&$\frac{1}{1200}$ &$\frac{-1}{1200}$ &$\frac{-7}{400}$&0&$\frac{7}{400}$&$\frac{1}{300}$&0&$\frac{-4}{75}$&$\frac{-7}{100}$&$\frac{28}{25}$&$\frac{1}{20}$&$\frac{21}{20}$&$\frac{2}{15}$&$\frac{-8}{3}$&$\frac{-112}{3}$\\
&&&&&&&&&&&&&&&\\
$0044$&$\frac{1}{2400}$ &$\frac{1}{800}$ &$\frac{3}{800}$&$\frac{-17}{600}$&$\frac{9}{800}$&$\frac{1}{50}$&$\frac{-51}{200}$&$\frac{8}{75}$&$\frac{9}{50}$&$\frac{24}{25}$&0&0&$\frac{2}{5}$&$\frac{28}{5}$&$\frac{112}{5}$\\
&&&&&&&&&&&&&&&\\
$0026$&$\frac{7}{19200}$ &$\frac{-7}{19200}$ &$\frac{-9}{6400}$&0&$\frac{9}{6400}$&$\frac{7}{300}$&0&$\frac{-28}{75}$&$\frac{-9}{100}$&$\frac{36}{25}$&$\frac{-7}{80}$&$\frac{-27}{80}$&$\frac{-19}{10}$&-7&-8\\
&&&&&&&&&&&&&&&\\
\hline
\end{tabular}
}
\end{center}
\vfill

\newpage 

\begin{center}
{\tiny
\textbf{Table 3.} (Theorem 2.2 (i) (contd.))\\
\begin{tabular}{|c|c|c|c|c|c|c|c|c|c|c|c|c|c|c|c|}
\hline
&&&&&&&&&&&&&&&\\
$ijkl$& $a_{16}$&$a_{17}$& $a_{18}$&$a_{19}$ &$a_{20}$ & $a_{21}$& $a_{22}$&$a_{23}$ &$a_{24}$& $a_{25}$& $a_{26}$&$a_{27}$ &$a_{28}$& $a_{29}$& $a_{30}$\\
\hline
&&&&&&&&&&&&&&&\\
$0062$&$\frac{-2432}{15}$&$\frac{-3}{2}$&-9&$\frac{-9}{2}$&-27&$\frac{2}{3}$&$\frac{10}{3}$&$\frac{68}{3}$&$\frac{-9}{4}$&$\frac{27}{4}$&$\frac{1}{2}$&7&$\frac{38}{15}$&$\frac{17}{12}$&$\frac{-7}{4}$\\
&&&&&&&&&&&&&&&\\
$0044$&$\frac{128}{5}$&$\frac{3}{2}$&6&$\frac{27}{2}$&54&$\frac{-1}{2}$& -6& -8&1&9&$\frac{-3}{2}$&-6&-2&0&1\\
&&&&&&&&&&&&&&&\\
$0026$&$\frac{32}{5}$&$\frac{-9}{16}$&$\frac{-3}{2}$&$\frac{-243}{16}$&$\frac{-81}{2}$&$\frac{17}{16}$&$\frac{5}{2}$& 8&$\frac{9}{16}$&$\frac{-243}{16}$&$\frac{21}{16}$&$\frac{3}{2}$&$\frac{19}{10}$&$\frac{-17}{16}$&$\frac{-21}{16}$\\
\hline
\end{tabular}
}
\end{center}

\newpage 

\begin{center}
{\tiny 

\textbf{Table 4.} (Theorem 2.2 (ii)) \\
$~~$ \\
\begin{tabular}{|c|c|c|c|c|c|c|c|c|c|c|c|c|c|c|}
\hline
&&&&&&&&&&&&&&\\
$ijkl$&$ b_1$ &$ b_2$& $ b_3$& $ b_4$&$ b_5$ &$ b_6$& $ b_7$& $ b_8$&$ b_9$ &$ b_{10}$ & $ b_{11}$&$ b_{12}$& $ b_{13}$&$ b_{14}$\\
\hline 
&&&&&&&&&&&&&&\\
$0521$& 0&$\frac{-26}{451}$& 0&$\frac{108}{451}$& $\frac{208}{451}$&0 &$\frac{864}{451}$&0&$\frac{-3112}{451}$&$\frac{16768}{451}$ 
&$\frac{324}{451}$&$\frac{98496}{451}$&$\frac{3504}{451}$& $\frac{832}{451}$\\
&&&&&&&&&&&&&&\\
$0341$& 0&$\frac{10}{451}$& 0&$\frac{72}{451}$& $\frac{160}{451}$&0 &$\frac{-1152}{451}$&0&$\frac{14744}{1353}$&$\frac{-32392}{451}$ 
&$\frac{-432}{451}$&$\frac{-64224}{451}$&$\frac{-7296}{451}$& $\frac{-320}{451}$\\
&&&&&&&&&&&&&&\\
$0323$ & 0&$\frac{-26}{451}$& 0&$\frac{108}{451}$& $\frac{104}{451}$&0 &$\frac{432}{451}$&0&$\frac{-4713}{451}$&$\frac{16768}{451}$ 
&$\frac{162}{451}$&$\frac{98496}{451}$&$\frac{5360}{451}$& $\frac{832}{451}$\\
\hline
\end{tabular}
}
\end{center}

\newpage 
\begin{center}
{\tiny
\textbf{Table 5.}(Theorem 2.2 (iii)) \\
$~~$\\
\begin{tabular}{|c|c|c|c|c|c|c|c|c|c|c|c|c|c|c|c|}
\hline
&&&&&&&&&&&&&&&\\
$ijkl$&$c_1$ &$c_2$& $c_3$& $c_4$&$c_5$ &$c_6$& $c_7$& $c_8$&$c_9$ &$c_{10}$& $c_{11}$& $c_{12}$&$c_{13}$ &$c_{14}$& $c_{15}$ \\
\hline 
&&&&&&&&&&&&&&&\\
$0071$&$\frac{1}{46}$&$\frac{4}{23}$&$\frac{4}{23}$&$\frac{1}{46}$&$\frac{-1}{46}$&$\frac{-4}{23}$&$\frac{4}{23}$&$\frac{1}{46}$&$\frac{1}{23}$&$\frac{64}{23}$&$\frac{64}{23}$&$\frac{1}{23}$&$\frac{14}{23}$&0&$\frac{28}{23}$\\
&&&&&&&&&&&&&&&\\
$0053$&$\frac{-1}{92}$&$\frac{3}{23}$&$\frac{-1}{23}$&$\frac{3}{92}$&$\frac{1}{92}$&$\frac{6}{23}$&$\frac{2}{23}$&$\frac{3}{92}$&$\frac{1}{23}$&$\frac{-96}{23}$&$\frac{32}{23}$&$\frac{-3}{23}$&$\frac{-5}{46}$&$\frac{38}{23}$&$\frac{191}{23}$\\
&&&&&&&&&&&&&&&\\
$0035$&$\frac{-1}{184}$&$\frac{9}{92}$&$\frac{1}{92}$&$\frac{-9}{184}$&$\frac{1}{184}$&$\frac{9}{23}$&$\frac{-1}{23}$&$\frac{-9}{184}$&$\frac{1}{23}$&$\frac{-144}{23}$&$\frac{-16}{23}$&$\frac{9}{23}$&$\frac{-5}{92}$&$\frac{98}{23}$&$\frac{-677}{92}$\\
\hline
\end{tabular}
}
\end{center}

\newpage 
\begin{center}
{\tiny
\textbf{Table 5.}(Theorem 2.2 (iii)) (contd.) \\
$~~$\\
\begin{tabular}{|c|c|c|c|c|c|c|c|c|c|c|c|c|c|c|c|}
\hline
&&&&&&&&&&&&&&&\\
$ijkl$& $c_{16}$&$c_{17}$& $c_{18}$&$c_{19}$ &$c_{20}$ & $c_{21}$& $c_{22}$&$c_{23}$ &$c_{24}$& $c_{25}$& $c_{26}$&$c_{27}$ &$c_{28}$& $c_{29}$& $c_{30}$\\
\hline
&&&&&&&&&&&&&&&\\
$0071$&$\frac{28}{23}$&0&$\frac{364}{23}$&$\frac{1092}{23}$&0&$\frac{-728}{23}$&0&$\frac{5824}{23}$&$\frac{5824}{23}$&-5&$\frac{7}{3}$&4&$\frac{28}{3}$&$\frac{80}{3}$&$\frac{-112}{3}$\\
&&&&&&&&&&&&&&&\\
$0053$&$\frac{-43}{23}$&$\frac{-78}{23}$&$\frac{-342}{23}$&$\frac{186}{23}$&$\frac{1224}{23}$&$\frac{204}{23}$&$\frac{1296}{23}$&$\frac{-5760}{23}$&$\frac{-4368}{23}$&$\frac{15}{4}$&$\frac{1}{4}$&$\frac{-15}{4}$&$\frac{-17}{4}$&-20& 12\\
&&&&&&&&&&&&&&&\\
$0035$&$\frac{595}{92}$&$\frac{-122}{23}$&$\frac{350}{23}$&$\frac{-776}{23}$&$\frac{-1136}{23}$&$\frac{380}{23}$&$\frac{-2464}{23}$&$\frac{3992}{23}$&$\frac{1720}{23}$&$\frac{-5}{8}$&$\frac{-21}{8}$&$\frac{5}{8}$&$\frac{9}{8}$&-2& 18\\
\hline
\end{tabular}
}
\end{center}
\newpage 
\begin{center}
{\tiny 
\textbf{Table 6. }(Theorem 2.2 (iv))
\begin{tabular}{|c|c|c|c|c|c|c|c|c|c|c|c|c|c|c|}
\hline
$ijkl$&$ d_1$ &$ d_2$& $ d_3$& $ d_4$&$ d_5$ &$ d_6$& $ d_7$& $ d_8$&$ d_9$ &$ d_{10}$\\
\hline 
$0512$& 0&$\frac{1}{261}$& $\frac{4}{261}$& 0& $\frac{12}{29}$&0 &0 &$\frac{-3}{29}$&$\frac{3234}{261}$&$\frac{16192}{261}$\\
&&&&&&&&&&\\
$0332$& 0&$\frac{1}{261}$& $\frac{-8}{261}$& 0& $\frac{8}{29}$&0 &0 &$\frac{1}{29}$&$\frac{-3788}{261}$&$\frac{-17648}{261}$\\
&&&&&&&&&&\\
$0314$& 0&$\frac{1}{261}$& $\frac{2}{261}$& 0& $\frac{6}{29}$&0 &0 &$\frac{-3}{29}$&$\frac{3800}{261}$&$\frac{16192}{261}$\\
&&&&&&&&&&\\
$0152$& 0&$\frac{1}{261}$& $\frac{16}{261}$& 0& $\frac{16}{87}$&0 &0 &$\frac{-1}{87}$&$\frac{1960}{261}$&$\frac{28432}{261}$\\
\end{tabular}

$ijkl$& $ d_{11}$&$ d_{12}$& $ d_{13}$&$ d_{14}$&  $ d_{15}$ & $ d_{16}$& $ d_{17}$&$ d_{18}$& $ d_{19}$&$ d_{20}$\\
\hline
&&&&&&&&&&\\
$0512$&$\frac{1258}{87}$& $\frac{5224}{87}$& $\frac{-9470}{261}$&$\frac{-38912}{261}$ &$\frac{-8264}{87}$&$\frac{-45728}{87}$&$\frac{-4844}{87}$& $\frac{-1184}{87}$& $\frac{57688}{261}$&$\frac{223456}{261}$\\
&&&&&&&&&&\\
$0332$ &$\frac{-1100}{87}$& $\frac{-6752}{87}$& $\frac{7078}{261}$&$\frac{1208}{261}$ &$\frac{10696}{87}$&$\frac{1808}{3}$&$\frac{2524}{87}$& $\frac{6880}{87}$& $\frac{-45608}{261}$&$\frac{-295952}{261}$\\
&&&&&&&&&&\\
$0314$ &$\frac{1151}{87}$& $\frac{5224}{87}$& $\frac{-8650}{261}$&$\frac{-38912}{261}$&$\frac{-10048}{87}$&$\frac{-45728}{87}$&$\frac{-3988}{87}$& $\frac{-1184}{87}$& $\frac{54944}{261}$&$\frac{223456}{261}$\\
&&&&&&&&&&\\
\hline
\end{tabular}
}
\end{center}

\newpage
\begin{center}
{\tiny
\textbf{Table 6.} (Theorem 2.2 (iv))(contd.) \\
$~~$ \\
\begin{tabular}{|c|c|c|c|c|c|c|c|c|}
\hline
&&&&&&&&\\
$ijkl$&$ d_{21}$ &$ d_{22}$& $ d_{23}$& $ d_{24}$&$ d_{25}$ &$ d_{26}$&$ d_{27}$& $ d_{28}$ \\
\hline
&&&&&&&&\\
$0512$&$\frac{26914}{261}$&$\frac{52480}{261}$&$\frac{14720}{261}$& $\frac{-52480}{261}$& $\frac{-3536}{261}$&$\frac{-14600}{261}$&$\frac{-3562}{261}$&$\frac{-14656}{261}$\\
&&&&&&&&\\
$0332$&$\frac{-21050}{261}$&$\frac{-45080}{261}$&$\frac{-20800}{261}$& $\frac{129536}{261}$& $\frac{3724}{261}$&$\frac{18160}{261}$&$\frac{3704}{261}$&$\frac{18176}{261}$\\
&&&&&&&&\\
$0314$&$\frac{26246}{261}$&$\frac{52480}{261}$&$\frac{15712}{261}$& $\frac{-52480}{261}$& $\frac{-3856}{261}$&$\frac{-14600}{261}$&$\frac{-3869}{261}$&$\frac{-14656}{261}$\\
\hline
\end{tabular}
}
\end{center}

\pagebreak
\begin{center}
{Table 7.}(Theorem 2.3 (i))
\hspace{0.5 cm}
{\tiny
\begin{tabular}{|c|c|c|c|c|c|c|c|c|c|c|}
\hline
&&&&&&&&&&\\
$ijklm$&$a'_1$ &$a'_2$& $a'_3$& $a'_4$&$a'_5$ &$a'_6$& $a'_7$& $a'_8$&$a'_9$ &$a'_{10}$\\
\hline 
&&&&&&&&&&\\
$00026$&0&$\frac{1}{600}$&0&$\frac{-1}{600}$&$\frac{-7}{200}$&$\frac{1}{300}$&$\frac{7}{200}$&$\frac{-4}{75}$&$\frac{-7}{100}$&$\frac{28}{25}$\\
&&&&&&&&&&\\ 
$00062$&0&$\frac{7}{2400}$&0&$\frac{-7}{2400}$&$\frac{-9}{800}$&$\frac{7}{300}$&$\frac{9}{800}$&$\frac{-28}{75}$&$\frac{-9}{100}$&$\frac{36}{25}$\\
&&&&&&&&&&\\ 
$00044$&0&$\frac{1}{600}$&0&$\frac{-1}{600}$&$\frac{3}{200}$&$\frac{-1}{150}$&$\frac{-3}{200}$&$\frac{8}{75}$&$\frac{-3}{50}$&$\frac{24}{25}$\\
\hline
\end{tabular}
}
\end{center}

\pagebreak

\newpage
\begin{center}
Table 7 (contd.)
\hspace{0.5 cm}
{\tiny
\begin{tabular}{|c|c|c|c|c|c|c|c|c|c|c|}
\hline
&&&&&&&&&&\\
$ijklm$& $a'_{11}$&$a'_{12}$& $a'_{13}$&$a'_{14}$ &$a'_{15}$ & $a'_{16}$& $a'_{17}$&$a'_{18}$& $a'_{19}$&$a'_{20}$\\
\hline
&&&&&&&&&&\\
$00026$&0&0&0&$\frac{4}{15}$&$\frac{-16}{3}$&$\frac{-512}{15}$&0&-3&0&-9\\
&&&&&&&&&&\\
$00062$&0&0&0&$\frac{-16}{5}$&-8&$\frac{32}{5}$&0&$\frac{-3}{2}$&0&$\frac{-81}{2}$\\
&&&&&&&&&&\\
$00044$&0&0&0&$\frac{8}{5}$&16&$\frac{128}{5}$&0&2&0&18\\
\hline
\end{tabular}
}
\end{center}

\begin{center}

Table 7. (contd.)
\hspace{0.5 cm}
{\tiny 
\begin{tabular}{|c|c|c|c|c|c|c|c|c|c|c|}
\hline
&&&&&&&&&&\\
$ijklm$& $a'_{21}$&$a'_{22}$& $a'_{23}$&$a'_{24}$ &$a'_{25}$ & $a'_{26}$& $a'_{27}$&$a'_{28}$& $a'_{29}$&$a'_{30}$\\
\hline
&&&&&&&&&&\\
$00026$&0&$\frac{4}{3}$&$\frac{20}{3}$&0&0&0&1&0&0&0\\
&&&&&&&&&&\\
$00062$&0&$\frac{5}{2}$&8&0&0&0&$\frac{3}{2}$&0&0&0\\
&&&&&&&&&&\\
$00044$&0&-2&-8&0&0&0&-2&0&0&0\\
\hline
\end{tabular}
}
\end{center}

\begin{center}

Table 7.(contd.)
\hspace{0.5 cm}
{\tiny 
\begin{tabular}{|c|c|c|c|c|c|c|c|c|c|c|}
\hline
&&&&&&&&&&\\
$ijklm$&$a'_1$ &$a'_2$& $a'_3$& $a'_4$&$a'_5$ &$a'_6$& $a'_7$& $a'_8$&$a'_9$ &$a'_{10}$\\
\hline 
&&&&&&&&&&\\
$20042$&$\frac{7}{9600}$&$\frac{-7}{9600}$&$\frac{-9}{3200}$&0&$\frac{9}{3200}$&$\frac{7}{300}$&0&$\frac{-28}{75}$&$\frac{-9}{100}$&$\frac{36}{25}$\\
&&&&&&&&&&\\
$20024$&$\frac{1}{2400}$&$\frac{-1}{2400}$&$\frac{3}{800}$&0&$\frac{-3}{800}$&$\frac{-1}{150}$&0&$\frac{8}{75}$&$\frac{-3}{50}$&$\frac{24}{25}$\\
&&&&&&&&&&\\
$31121$&$\frac{13}{4800}$&$\frac{-13}{4800}$&$\frac{9}{1600}$&0&$\frac{-9}{1600}$&$\frac{-13}{600}$&0&$\frac{26}{75}$&$\frac{-9}{200}$&$\frac{18}{25}$\\
\hline
\end{tabular}
}
\end{center}

\begin{center}
Table 7.(contd.)
\hspace{0.5 cm}
{\tiny 
\begin{tabular}{|c|c|c|c|c|c|c|c|c|c|c|}
\hline
&&&&&&&&&&\\
$ijklm$& $a'_{21}$&$a'_{22}$& $a'_{23}$&$a'_{24}$ &$a'_{25}$ & $a'_{26}$& $a'_{27}$&$a'_{28}$& $a'_{29}$&$a'_{30}$\\
\hline
&&&&&&&&&&\\
$20042$&$\frac{-3}{8}$&3&4&0&0&$\frac{9}{8}$&3&1&1&0\\
&&&&&&&&&&\\
$20024$&$\frac{1}{2}$&0&-8&1&-9&$\frac{1}{2}$&-2&0&0&1\\
&&&&&&&&&&\\
$31121$&$\frac{-1}{4}$&1&-14&2&0&0&$\frac{-3}{2}$&-1&$\frac{1}{2}$&$\frac{3}{2}$\\
&&&&&&&&&&\\
$30212$&0&4&-8&$\frac{3}{2}$&$\frac{-9}{2}$&1&-4&-1&0,&$\frac{3}{2}$\\
&&&&&&&&&&\\
$32012$&$\frac{1}{2}$&6&4&0&0&$\frac{3}{2}$&3&1&1&0\\
&&&&&&&&&&\\
$30032$&$\frac{-1}{4}$&4&-4&0&0&$\frac{3}{2}$&6&$\frac{3}{2}$&$\frac{3}{2}$&0\\
&&&&&&&&&&\\
$30014$&1&4&-8&1&-9&1&-2&0&0&1\\\
&&&&&&&&&&\\
$41111$&0&5&-14&2&0&$\frac{3}{4}$&$\frac{-3}{2}$&-1&$\frac{1}{2}$&$\frac{3}{2}$\\
&&&&&&&&&&\\
$40022$&$\frac{1}{2}$&6&-12&0&0&$\frac{3}{2}$&9&2&2&0\\
&&&&&&&&&&\\
$50012$&2&10&-12&0&0&$\frac{3}{2}$&9&2&2&0\\
&&&&&&&&&&\\
\hline
\end{tabular}
}
\end{center}

\begin{center}
{\tiny
Table 8.(Theorem 2.3(ii))
\hspace{0.5 cm}
\begin{tabular}{|c|c|c|c|c|c|c|c|c|c|c|}
\hline 
&&&&&&&&&&\\
$ijklm$&$b'_1$ &$b'_2$& $b'_3$& $b'_4$&$b'_5$ &$b'_6$& $b'_7$& $b'_8$&$b'_9$ &$b'_{10}$\\
\hline
&&&&&&&&&&\\
$00323$&0 &$\frac{4}{451}$&0&$\frac{78}{451}$&$\frac{32}{451}$&0&$\frac{-624}{451}$&0&$\frac{44746}{4059}$&$\frac{-17944}{4059}$\\
&&&&&&&&&&\\
$00143$&0 &$\frac{-8}{451}$&0&$\frac{90}{451}$&$\frac{16}{451}$&0&$\frac{180}{451}$&0&$\frac{-8177}{1353}$&$\frac{3544}{1353}$\\
&&&&&&&&&&\\
$00341$&0 &$\frac{-8}{451}$&0&$\frac{90}{451}$&$\frac{32}{451}$&0&$\frac{360}{451}$&0&$\frac{-13648}{1353}$&$\frac{5992}{451}$\\
\hline
\end{tabular}
}
\end{center}

\begin{center}
{\tiny
Table 8.(contd.)
\hspace{0.5 cm}
\begin{tabular}{|c|c|c|c|c|c|c|c|c|c|c|}
\hline
&&&&&&&&&&\\
$ijklm$& $b'_{11}$&$b'_{12}$& $b'_{13}$&$b'_{14}$ &$b'_{15}$ & $b'_{16}$& $b'_{17}$&$b'_{18}$& $b'_{19}$&$b'_{20}$\\
\hline
&&&&&&&&&&\\
$00323$&$\frac{-234}{451}$&$\frac{-15456}{451}$&$\frac{-7232}{451}$&$\frac{-128}{451}$&$\frac{312}{451}$&$\frac{-2496}{451}$&$\frac{-2582}{369}$&$\frac{1664}{369}$&$\frac{-3682}{123}$&$\frac{2776}{123}$\\
&&&&&&&&&&\\
$00143$&$\frac{-1962}{451}$&$\frac{17136}{451}$&$\frac{3600}{451}$&$\frac{256}{451}$&$\frac{8028}{451}$&$\frac{-2880}{451}$&$\frac{499}{123}$&$\frac{-344}{123}$&$\frac{589}{41}$&$\frac{-648}{41}$\\
&&&&&&&&&&\\
$00341$&$\frac{4194}{451}$&$\frac{49608}{451}$&$\frac{7200}{451}$&$\frac{256}{451}$&$\frac{-180}{451}$&$\frac{-2880}{451}$&$\frac{752}{123}$&$\frac{-552}{41}$&$\frac{932}{41}$&$\frac{-2616}{41}$\\
\hline
\end{tabular}
}
\end{center}

\begin{center}
{\tiny
Table 8.(contd.)
\hspace{0.5 cm}
\begin{tabular}{|c|c|c|c|c|c|c|c|c|}
\hline
&&&&&&&&\\
$ijklm$&$b'_{21}$ &$b'_{22}$& $b'_{23}$& $b'_{24}$&$b'_{25}$ &$b'_{26}$& $b'_{27}$& $b'_{28}$ \\
\hline
&&&&&&&&\\
$00323$&$\frac{650}{123}$&$\frac{-1520}{123}$&$\frac{23752}{369}$&$\frac{-13888}{369}$&$\frac{4488}{41}$&$\frac{-6976}{41}$&$\frac{-1154}{123}$&$\frac{1016}{123}$\\
&&&&&&&&\\
$00143$&$\frac{-83}{41}$&$\frac{248}{41}$&$\frac{-2984}{123}$&$\frac{3136}{123}$&$\frac{-2184}{41}$&$\frac{3456}{41}$&$\frac{163}{41}$&$\frac{-240}{41}$\\
&&&&&&&&\\
$00341$&$\frac{-84}{41}$&$\frac{904}{41}$&$\frac{-4984}{123}$&$\frac{4544}{41}$&$\frac{-3384}{41}$&$\frac{11328}{41}$&$\frac{244}{41}$&$\frac{-896}{41}$\\
\hline
\end{tabular}
}
\end{center}

\begin{center}
{\tiny
Table 8.(contd.)
\hspace{0.5 cm}
\begin{tabular}{|c|c|c|c|c|c|c|c|c|c|c|}
\hline
&&&&&&&&&&\\
$ijklm$&$b'_1$ &$b'_2$& $b'_3$& $b'_4$&$b'_5$ &$b'_6$& $b'_7$& $b'_8$&$b'_9$ &$b'_{10}$\\
\hline 
&&&&&&&&&&\\
$ 20141$ &0&$\frac{28}{451}$&0&$\frac{54}{451}$&$\frac{112}{451}$&0&$\frac{-216}{451}$&0&$\frac{3483}{902}$&$\frac{-4736}{451}$\\
&&&&&&&&&&\\
$ 20321$ &0&$\frac{-8}{451}$&0&$\frac{90}{451}$&$\frac{128}{451}$&0&$\frac{1440}{451}$&0&$\frac{-26630}{1353}$&$\frac{9600}{451}$\\
&&&&&&&&&&\\
$ 31022$ &0&$\frac{-26}{451}$&0&$\frac{108}{451}$&$\frac{208}{451}$&0&$\frac{864}{451}$&0&$\frac{-6269}{451}$&$\frac{-1272}{451}$\\
\hline
\end{tabular}
}
\end{center}
\smallskip
\smallskip
\smallskip
\smallskip
\smallskip
\smallskip
\smallskip
\smallskip
\begin{center}
{\tiny
Table 8.(contd.)
\hspace{0.5 cm}
\begin{tabular}{|c|c|c|c|c|c|c|c|c|}
\hline
&&&&&&&&\\
$ijklm$&$b'_{21}$ &$b'_{22}$& $b'_{23}$& $b'_{24}$&$b'_{25}$ &$b'_{26}$& $b'_{27}$& $b'_{28}$\\
\hline
&&&&&&&&\\
$ 20141$ &$\frac{-219}{82}$&$\frac{-1032}{41}$&$\frac{-348}{41}$&$\frac{-4096}{41}$&$\frac{-36}{41}$&$\frac{-12096}{41}$&$\frac{191}{82}$&$\frac{840}{41}$\\
&&&&&&&&\\
$ 20321$ &$\frac{-172}{41}$&$\frac{1560}{41}$&$\frac{-13376}{123}$&$\frac{12416}{41}$&$\frac{-7632}{41}$&$\frac{19200}{41}$&$\frac{648}{41}$&$\frac{-1880}{41}$\\
&&&&&&&&\\
$ 31022$ &$\frac{-177}{41}$&$\frac{-96}{41}$&$\frac{-3528}{41}$&$\frac{-320}{41}$&$\frac{-5760}{41}$&$\frac{-576}{41}$&$\frac{561}{41}$&$\frac{40}{41}$\\
&&&&&&&&\\
$ 30113$ &$\frac{-49}{41}$&$\frac{576}{41}$&$\frac{-8128}{123}$&$\frac{16256}{123}$&$\frac{-4680}{41}$&$\frac{11328}{41}$&$\frac{443}{41}$&$\frac{-896}{41}$\\
&&&&&&&&\\
$ 30311$ &$\frac{-262}{41}$&$\frac{1560}{41}$&$\frac{-5200}{41}$&$\frac{12416}{41}$&$\frac{-9360}{41}$&$\frac{19200}{41}$&$\frac{804}{41}$&$\frac{-1880}{41}$\\
&&&&&&&&\\
$ 32111$ &$\frac{177}{41}$&$\frac{-1032}{41}$&$\frac{1232}{41}$&$\frac{-4096}{41}$&$\frac{2808}{41}$&$\frac{-12096}{41}$&$\frac{-151}{41}$&$\frac{840}{41}$\\
&&&&&&&&\\
$ 30131$ &$\frac{-561}{82}$&$\frac{-1524}{41}$&$\frac{-1024}{41}$&$\frac{-10656}{41}$&$\frac{-1548}{41}$&$\frac{-12096}{41}$&$\frac{505}{82}$&$\frac{1496}{41}$\\
&&&&&&&&\\
$ 41012$ &$\frac{-108}{41}$&$\frac{-96}{41}$&$\frac{-3120}{41}$&$\frac{-320}{41}$&$\frac{-5616}{41}$&$\frac{-576}{41}$&$\frac{548}{41}$&$\frac{40}{41}$\\
&&&&&&&&\\
$ 40121$ &$\frac{-192}{41}$&$\frac{-2016}{41}$&$\frac{-80}{41}$&$\frac{-17216}{41}$&$\frac{-144}{41}$&$\frac{-12096}{41}$&$\frac{136}{41}$&$\frac{2152}{41}$\\
&&&&&&&&\\
$ 50111$ &$\frac{354}{41}$&$\frac{-2016}{41}$&$\frac{3120}{41}$&$\frac{-17216}{41}$&$\frac{5616}{41}$&$\frac{-12096}{41}$&$\frac{-384}{41}$&$\frac{2152}{41}$\\
&&&&&&&&\\
\hline
\end{tabular}
}
\end{center}

\begin{center}
{\tiny
Table 9.(Theorem 2.3(iii))
\hspace{0.5 cm}
\begin{tabular}{|c|c|c|c|c|c|c|c|c|c|c|}
\hline
&&&&&&&&&&\\
$ijklm$&$c'_1$ &$c'_2$& $c'_3$& $c'_4$&$c'_5$ &$c'_6$& $c'_7$& $c'_8$&$c'_9$ &$c'_{10}$\\
\hline 
&&&&&&&&&&\\
$00512$&0&$\frac{2}{23}$&$\frac{2}{23}$&0&0&0&0&0&$\frac{1}{23}$&0\\
&&&&&&&&&&\\ 
$00134$&0&$\frac{3}{92}$&$\frac{-1}{92}$&0&0&0&0&0&$\frac{1}{23}$&0\\
&&&&&&&&&&\\ 
$00332$&0&$\frac{3}{46}$&$\frac{-1}{46}$&0&0&0&0&0&$\frac{1}{23}$&0\\
\hline
\end{tabular}
}
\end{center}

\begin{center}
{\tiny
Table 9.(contd.)
\hspace{0.5 cm}
\begin{tabular}{|c|c|c|c|c|c|c|c|c|c|c|}
\hline
&&&&&&&&&&\\ 
$ijklm$& $c'_{11}$&$c'_{12}$& $c'_{13}$&$c'_{14}$ &$c'_{15}$ & $c'_{16}$& $c'_{17}$&$c'_{18}$& $c'_{19}$&$c'_{20}$\\
\hline
&&&&&&&&&&\\ 
$00512$&0&$\frac{1}{23}$&$\frac{19}{23}$&0&$\frac{22}{23}$&$\frac{22}{23}$&0&16&24&0\\
&&&&&&&&&&\\ 
$00134$&0&$\frac{-3}{23}$&$\frac{-1}{46}$&$\frac{-8}{23}$&$\frac{173}{46}$&$\frac{17}{46}$&0&0&0&12\\
&&&&&&&&&&\\ 
$00332$&0&$\frac{-3}{23}$&$\frac{-1}{23}$&$\frac{-16}{23}$&$\frac{277}{46}$&$\frac{-35}{46}$&0&-6&6&24\\
\hline
\end{tabular}
}
\end{center}

\begin{center}
{\tiny
Table 9.(contd.)
\hspace{0.5 cm}
\begin{tabular}{|c|c|c|c|c|c|c|c|c|c|c|}
\hline
&&&&&&&&&&\\ 
$ijklm$& $c'_{21}$&$c'_{22}$& $c'_{23}$&$c'_{24}$ &$c'_{25}$ & $c'_{26}$& $c'_{27}$&$c'_{28}$& $c'_{29}$&$c'_{30}$\\
\hline
&&&&&&&&&&\\ 
$00512$&$\frac{-416}{23}$&0&$\frac{2912}{23}$&$\frac{2912}{23}$&$\frac{-5}{2}$&$\frac{7}{6}$&$\frac{3}{2}$&$\frac{31}{6}$&$\frac{40}{3}$&$\frac{-56}{3}$\\
&&&&&&&&&& \\ 
$00134$ & $\frac{140}{23}$&$\frac{272}{23}$&$\frac{-1312}{23}$&$\frac{-1072}{23}$&$\frac{3}{4}$&$\frac{-1}{4}$&$\frac{-3}{4}$&$\frac{-7}{4}$&-4&4\\
&&&&&&&&&&\\ 
$00332$&$\frac{140}{23}$&$\frac{272}{23}$&$\frac{-2416}{23}$&$\frac{-2176}{23}$&$\frac{3}{2}$&0&$\frac{-3}{2}$&-3&-8&8\\
\hline
\end{tabular}
}
\end{center}

\begin{center}
{\tiny
Table 9.(contd.)
\hspace{0.5 cm}
\begin{tabular}{|c|c|c|c|c|c|c|c|c|c|c|}
\hline
&&&&&&&&&&\\ 
$ijklm$&$c'_1$ &$c'_2$& $c'_3$& $c'_4$&$c'_5$ &$c'_6$& $c'_7$& $c'_8$&$c'_9$ &$c'_{10}$\\
\hline 
&&&&&&&&&&\\ 
$21221$&0&$\frac{9}{23}$&$\frac{-1}{23}$&0&0&0&0&0&$\frac{1}{23}$&0\\
&&&&&&&&&&\\ 
$23021$&0&$\frac{27}{46}$&$\frac{1}{46}$&0&0&0&0&0&$\frac{1}{23}$&0\\
&&&&&&&&&&\\ 
$33011$&0&$\frac{27}{23}$&$\frac{1}{23}$&0&0&0&0&0&$\frac{1}{23}$&0\\
\hline
\end{tabular}
}
\end{center}

\begin{center}
{\tiny
Table 9.(contd.)
\hspace{0.5 cm}
\begin{tabular}{|c|c|c|c|c|c|c|c|c|c|c|}
\hline
&&&&&&&&&&\\ 
$ijklm$&$c'_{21}$ &$c'_{22}$& $c'_{23}$& $c'_{24}$&$c'_{25}$ &$c'_{26}$& $c'_{27}$& $c'_{28}$& $c'_{29}$& $c'_{30}$\\
\hline
&&&&&&&&&&\\ 
$21221$&$\frac{192}{23}$&$\frac{128}{23}$&$\frac{-2496}{23}$&$\frac{-2816}{23}$&3&1&-1&-3&-16&0\\
&&&&&&&&&&\\ 
$23021$&$\frac{64}{23}$&$\frac{-1216}{23}$&$\frac{3136}{23}$&$\frac{1024}{23}$&0&1&2&3&-4&-12\\
&&&&&&&&&&\\ 
$33011$&$\frac{64}{23}$&$\frac{-1216}{23}$&$\frac{3136}{23}$&$\frac{1024}{23}$&0&1&2&3&-4&-12\\
&&&&&&&&&&\\ 
$30122$&$\frac{-608}{23}$&$\frac{192}{23}$&$\frac{4000}{23}$&$\frac{2784}{23}$&$\frac{-1}{2}$&$\frac{1}{2}$&$\frac{7}{2}$&$\frac{9}{2}$&8&-24\\
&&&&&&&&&&\\ 
$31013$&$\frac{-176}{23}$&$\frac{864}{23}$&$\frac{-1760}{23}$&$\frac{-1344}{23}$&2&-1&0&-3&-12&12\\
&&&&&&&&&&\\ 
$31211$&$\frac{192}{23}$&$\frac{128}{23}$&$\frac{-2496}{23}$&$\frac{-2816}{23}$&3&1&-1&-3&-16&0\\
&&&&&&&&&&\\ 
$31031$&$\frac{-304}{23}$&$\frac{-848}{23}$&$\frac{4976}{23}$&$\frac{1760}{23}$&0&$\frac{3}{2}$&3&$\frac{9}{2}$&-6&-18\\
&&&&&&&&&&\\ 
$40112$&$\frac{-608}{23}$&$\frac{192}{23}$&$\frac{4000}{23}$&$\frac{2784}{23}$&$\frac{-1}{2}$&$\frac{1}{2}$&$\frac{7}{2}$&$\frac{9}{2}$&8&-24\\
&&&&&&&&&&\\ 
$41021$&$\frac{-672}{23}$&$\frac{-480}{23}$&$\frac{6816}{23}$&$\frac{2496}{23}$&0&2&4&6&-8&-24\\
&&&&&&&&&&\\ 
$51011$&$\frac{-672}{23}$&$\frac{-480}{23}$&$\frac{6816}{23}$&$\frac{2496}{23}$&0&2&4&6&-8&-24\\
&&&&&&&&&&\\ 
\hline
\end{tabular}
}
\end{center}

\begin{center}
{\tiny
Table 10.(Theorem 2.3(iv))
\hspace{0.5 cm}
\begin{tabular}{|c|c|c|c|c|c|c|c|c|c|c|}
\hline
&&&&&&&&&&\\
$ijklm$&$d'_1$ &$d'_2$& $d'_3$& $d'_4$&$d'_5$ &$d'_6$& $d'_7$& $d'_8$&$d'_9$ &$d'_{10}$\\
\hline 
&&&&&&&&&&\\
$00035$&0&$\frac{1}{261}$&0& $\frac{32}{261}$&0&$\frac{32}{87}$&0&$\frac{1}{87}$&0&$\frac{-376}{87}$\\
&&&&&&&&&&\\
$00071$&0&$\frac{1}{261}$&0&$\frac{8}{261}$&0&$\frac{24}{29}$&0&$\frac{3}{29}$&0&$\frac{168}{29}$\\
&&&&&&&&&&\\
$00053$&0&$\frac{1}{261}$&0&$\frac{-16}{261}$&0&$\frac{16}{29}$&0&$\frac{-1}{29}$&0&$\frac{-1376}{87}$\\
\hline
\end{tabular}
}
\end{center}
\begin{center}
{\tiny
Table 10.(contd.)
\hspace{0.5 cm}
\begin{tabular}{|c|c|c|c|c|c|c|c|c|c|c|}
\hline
&&&&&&&&&&\\
$ijklm$& $d'_{11}$&$d'_{12}$& $d'_{13}$&$d'_{14}$ &$d'_{15}$ & $d'_{16}$& $d'_{17}$&$d'_{18}$& $d'_{19}$&$d'_{20}$\\
\hline
&&&&&&&&&&\\
$00035$&0&$\frac{8}{29}$&0& $\frac{632}{87}$&0&$\frac{352}{29}$&0&$\frac{1328}{29}$&0&$\frac{-4192}{87}$\\
&&&&&&&&&&\\
$00071$&0&$\frac{420}{29}$&0&$\frac{-1064}{29}$&0&$\frac{-672}{29}$&0&$\frac{-1680}{29}$&0& $\frac{5600}{29}$\\
&&&&&&&&&&\\
$00053$&0&$\frac{-292}{29}$&0&$\frac{1984}{87}$&0&$\frac{3840}{29}$&0&$\frac{544}{29}$&0& $\frac{-16736}{87}$\\
\hline
\end{tabular}
}
\end{center}

\begin{center}
{\tiny
Table 10.(contd.)
\hspace{0.5 cm}
\begin{tabular}{|c|c|c|c|c|c|c|c|c|}
\hline
&&&&&&&&\\
$ijklm$&$d'_{21}$ &$d'_{22}$& $d'_{23}$& $d'_{24}$&$d'_{25}$ &$d'_{26}$& $d'_{27}$& $d'_{28}$ \\
\hline
&&&&&&&&\\
$00035$&0&$\frac{-9544}{87}$&0&$\frac{16384}{87}$&0&$\frac{332}{87}$&0&$\frac{328}{87}$\\
&&&&&&&&\\
$00071$&0&$\frac{2744}{29}$&0&$\frac{1792}{29}$&0&$\frac{-196}{29}$&0& $\frac{-196}{29}$\\
&&&&&&&&\\
$00053$&0&$\frac{-1664}{87}$&0&$\frac{-12160}{87}$&0&$\frac{1336}{87}$&0& $\frac{1316}{87}$\\
\hline
\end{tabular}
}
\end{center}

\begin{center}
{\tiny
Table 10.(contd.)
\hspace{0.5 cm}
\begin{tabular}{|c|c|c|c|c|c|c|c|c|c|c|}
\hline
$ijklm$&$d'_1$ &$d'_2$& $d'_3$& $d'_4$&$d'_5$ &$d'_6$& $d'_7$& $d'_8$&$d'_9$ &$d'_{10}$\\
\hline 
&&&&&&&&&&\\
$22031$&0&$\frac{1}{261}$&$\frac{-4}{261}$&0&$\frac{12}{29}$&0&0&$\frac{3}{29}$&$\frac{2786}{261}$&$\frac{-328}{9}$\\
&&&&&&&&&&\\
$20213$&0&$\frac{1}{261}$&$\frac{-16}{261}$&0&$\frac{16}{87}$&0&0&$\frac{1}{87}$&$\frac{-592}{261}$&$\frac{-19664}{261}$\\
&&&&&&&&&&\\
$20033$&0&$\frac{1}{261}$&$\frac{4}{261}$&0&$\frac{4}{29}$&0&0&$\frac{-1}{29}$&$\frac{-734}{261}$&$\frac{11800}{261}$\\
&&&&&&&&&&\\
\hline
\end{tabular}
}
\end{center}

\begin{center}
{\tiny
Table 10.(contd.)
\hspace{0.5 cm}
\begin{tabular}{|c|c|c|c|c|c|c|c|c|}
\hline
&&&&&&&&\\
$ijklm$&$d'_{21}$ &$d'_{22}$& $d'_{23}$& $d'_{24}$&$d'_{25}$ &$d'_{26}$& $d'_{27}$& $d'_{28}$ \\
\hline
&&&&&&&&\\
$22031$&$\frac{12416}{261}$&$\frac{-44792}{261}$& $\frac{-12416}{261}$&$\frac{86912}{261}$&$\frac{-1846}{261}$&$\frac{11572}{261}$&$\frac{-1874}{261}$&$\frac{9536}{261}$\\ 
&&&&&&&&\\
$20213$&$\frac{-8308}{261}$&$\frac{-109208}{261}$& $\frac{64}{261}$&$\frac{132608}{261}$&$\frac{1604}{261}$&$\frac{18616}{261}$&$\frac{1588}{261}$&$\frac{17576}{261}$\\ 
&&&&&&&&\\
$20033$&$\frac{-304}{9}$&$\frac{52768}{261}$&$\frac{8384}{261}$&$\frac{-95104}{261}$&$\frac{1738}{261}$&$\frac{-10748}{261}$&$\frac{1730}{261}$&$\frac{-12856}{261}$\\
\hline
\end{tabular}
}
\end{center}

\end{document}